\documentclass{article}
\usepackage{amsmath}[1996/11/01]
\usepackage{amssymb,amsthm,amsxtra}

\newtheorem{thm}{Theorem}
\numberwithin{thm}{section}
\newtheorem{cor}[thm]{Corollary}
\newtheorem{lem}[thm]{Lemma}
\newtheorem{prop}[thm]{Proposition}

\theoremstyle{definition}

\theoremstyle{remark}
\newtheorem*{rem}{Remark}
\newtheorem{rems}{Remark}[thm]

\newcommand{\Z}{\mathbb Z}
\newcommand{\C}{\mathbb C}

\newcommand{\R}{\mathbb R}

\newcommand{\Gammae}{\Gamma_{\!e}}
\newcommand{\Gammat}{\Gamma_{\!t}}
\newcommand{\Gammar}{\Gamma_{\!r}}
\newcommand{\Gammah}{\Gamma_{\!h}}
\DeclareMathOperator\Res{Res}
\DeclareMathOperator\sgn{sgn}
\newcommand{\II}{\mathord{I\!I}}
\newcommand{\blambda}{{\boldsymbol\lambda}}
\newcommand{\bmu}{{\boldsymbol\mu}}
\newcommand{\PVint}{\text{\rm PV}\!\!\int}

\begin{document}

\title{Limits of elliptic hypergeometric integrals}
\author{Eric M. Rains}
\date{September 5, 2007}
\maketitle
\begin{abstract}
In \cite{Rains:Transformations}, the author proved a number of multivariate
elliptic hypergeometric integrals.  The purpose of the present note is to
explore more carefully the various limiting cases (hyperbolic,
trigonometric, rational, and classical) that exist.  In particular, we show
(using some new estimates of generalized gamma functions) that the
hyperbolic integrals (previously treated as purely formal limits) are
indeed limiting cases.  We also obtain a number of new trigonometric
($q$-hypergeometric) integral identities as limits from the elliptic level.
\end{abstract}

\tableofcontents

\section{Introduction}

In \cite{Rains:Transformations}, the author proved and extended a
pair of multivariate elliptic hypergeometric integrals conjectured
by van Diejen and Spiridonov \cite{vanDiejenJF/SpiridonovVP:2001},
elliptic analogues of integrals due to Gustafson \cite{GustafsonRA:1992}.
In a follow-up paper \cite{vanDiejenJF/SpiridonovVP:2005}, van
Diejen and Spiridonov proved ``hyperbolic'' degenerations of these
integrals (after the hyperbolic gamma function of Ruijsensaars
\cite{RuijsenaarsSNM:1999}, and results of Stokman \cite {StokmanJV:2005}
at the univariate level).  Unfortunately, the asymptotic estimates
they had available were insufficient to derive these hyperbolic
integrals as limits of the elliptic integrals; instead, they were
forced to degenerate the elliptic proofs.  The objective of the
present paper is to put these (and other) degenerations on a sounder
footing by showing that they are indeed limiting cases of the
elliptic integral.  It is to be hoped that a better understanding
of the relation between the elliptic and other integrals will lead
to new results at various levels, both by clarifying how arguments
at low levels extend to higher levels and by providing new results
at low levels as limits of the elliptic identities.  The latter
hope has already been fulfilled to some extent; see, for instance,
Corollary \ref{cor:triglim_IC12} below.

A full study of degenerations of the elliptic hypergeometric integrals
is beyond the scope of the present paper.  Here, we focus only on
the top level of each of the five types of limit (hyperbolic,
trigonometric, elliptic, rational, and classical); degenerations
within each type will be deferred to future work with F. van de
Bult.  We also do not consider in any detail the degenerations
between types; rather, we obtain each case directly as a limit from
the elliptic level.  (Note, however, that in several cases our
estimates are sufficiently uniform that one could view certain
between-level limits as special cases of the limits from the elliptic
level.) Similarly, we do not consider the discrete degenerations,
a.k.a. hypergeometric sums; aside from finite sums (which have
already been considered), the only interesting discrete degenerations
appear to arise via lower-level (e.g., hyperbolic or trigonometric)
integrals.  Finally, we comment that at the univariate level, there
are important relations between the trigonometric and hyperbolic
integrals (\cite{StokmanJV:2005}, also
\cite{vdBultF/RainsEM/StokmanJV:2007}); presumably, these have
analogues at the multivariate level, but the existing techniques
do not appear to apply.  (However, see \cite{StokmanJV:2007} for recent
progress on a slightly degenerate, but multivariate case.)

The plan of the paper is as follows.  In Section 2, we define the relevant
gamma functions (elliptic, hyperbolic, trigonometric, and rational), and
prove a number of asymptotic results relating them.  In order to make the
derivation of limit integrals as simple as possible, we have tried to make
these estimates as uniform as possible, and the relevant error estimates
similarly strong; since our techniques also apply equally well to
higher-order multiple gamma functions, we prove our estimates in that level
of generality.  Even for the ordinary (i.e., $r=2$) elliptic and hyperbolic
gamma functions, our results are new, as the results in the literature were
nonuniform (compare, for instance, Theorem \ref{thm:hyper_to_rat} to its
special cases Corollaries \ref{cor:hyper_to_rat} and
\ref{cor:hyper_to_rat2}, which were proved (for $r=2$, and not uniform) in
\cite{RuijsenaarsSNM:1999}).  Corollary \ref{cor:trig_to_rat}, a uniform
version of the $q\to 1$ asymptotics of $q$-symbols, may also be of
independent interest.

Section 3 considers the inequalities required to pinpoint where the
integrands of interest are maximized.  It turns out that these all follow
from a single master inequality (Lemma \ref{lem:gen_tri_ell}), which in
turn follows from an asymptotic analysis of an elliptic analogue of the
Cauchy determinant.  The result is an inequality stating that certain
combinations of elliptic gamma functions are exponentially small unless
their arguments alternate around the unit circle.  It may be worth
investigating other elliptic determinant and pfaffian identities to see
whether they give rise to interesting inequalities.

Section 4 begins the study of limit integrals with the hyperbolic case.
This case turns out to be relatively straightforward, given the estimates
and inequalities already established; in each case, a standard
tail-exchange argument gives the desired limit.  In fact, not only do we
obtain the hyperbolic integrals as limits of the corresponding elliptic
integrals, but moreover obtain exponentially small error estimates.  We
also discuss the extra arguments required to include the case in which the
integrands are multiplied by appropriate abelian functions (the
interpolation and biorthogonal functions of \cite{Rains:Transformations}).

Section 5 considers the trigonometric (i.e., basic hypergeometric) case.
It turns out that in addition to the integrals of Gustafson that originally
motivated van Diejen and Spiridonov, there are additional limiting cases.
That these should exist is suggested by the fact that in the
transformations of \cite{Rains:Transformations}, only one side typically
has a straightforward limit as $p\to 0$.  It turns out, however, that if
one breaks the symmetry of the integrand in a suitable way, one can arrange
for both sides of the transformations to have reasonable limits.  As a
special case (violating our rule of considering only the top level of each
case), we find that not only are the ordinary Macdonald polynomials limits
of the biorthogonal abelian functions, but in such a way that their
orthogonality follows as well; in particular, the Macdonald ``conjectures''
(proved in \cite{MacdonaldIG:1995}) for ordinary Macdonald polynomials are
limiting cases of the corresponding identities for biorthogonal abelian
functions.

The remaining two cases correspond to the case $q\to 1$ with $p$ fixed.  In
the rational case, considered in Section 6, the parameters all tend to 1,
and one obtains a hybrid of the hyperbolic and trigonometric cases.  In
particular, to obtain the full spectrum of possibilities, one must first
break the symmetry of the integrand, and then do a tail-exchange asymptotic
argument.  The arguments are still fairly straightforward, although the
resulting estimates are much weaker.

Finally, in the classical case, considered in Section 7, the parameters
behave in such a way that the integrand is exponentially small unless
certain inequalities are satisfied.  This gives rise to multivariate
analogues of the ordinary beta integral; notably those of Dixon
\cite{DixonAL:1905} and Selberg \cite{SelbergA:1944}.  Curiously, the most
natural forms of these limits remain elliptic in nature, with the integrand
involving powers of theta functions, although these can be removed with a
suitable change of variables.

The author would like to thank P. Forrester, J. Stokman and F. van
de Bult for motivating conversations regarding the trigonometric
and hyperbolic cases, and R. Askey for suggesting the use of the
modular transformation to derive classical limits (as in
\cite{NasrallahB/RahmanM:1985}), which led the author to consider
the paper \cite{NarukawaA:2004}; the author would also like to thank
an anonymous referee for pointing out that the original version of
Theorem \ref{thm:biorth_lim} was badly stated.  The author was supported
in part by NSF Grant No. DMS-0401387.

\subsubsection*{Conventions}
We use standard conventions for big $O$ notation in uniform estimates; that
is, we say that $f(v,z)=O(g(v,z))$ as $v\to 0^+$ uniformly over the region
$z\in D(v)$ (where $D(v)$ is a set depending on $v$) if there exist
constants $\delta>0$ and $C>0$ such that
\[
|f(v,z)|<C |g(v,z)|
\]
whenever $0<v<\delta$, $z\in D(v)$; and similarly for limits $v\to
\pm\infty$ (which can be viewed as limits $\pm 1/v\to 0^+$).

All logarithms are taken on the principal branch, with branch cut along the
negative real axis.  All powers are determined correspondingly; in
particular, square roots are chosen with positive real part.  (The obvious
exception is $\sqrt{-1}$, which is always taken to have positive imaginary
part; we refrain from denoting it by $i$ so as to retain $i$ as an indexing
variable.)

Finally, for a real number $x$, $\{x\}$ denotes the fractional part of $x$;
i.e., the unique representative in $[0,1)$ of $x+\Z$.

\section{Asymptotics of multiple gamma functions}
For integers $r\ge 0$, let $\phi^{(r)}(w;x;\omega_1,\dots,\omega_r)$ denote
the function
\[
\phi^{(r)}(w;x;\omega_1,\dots,\omega_r)
:=
\frac{e(xw)}{\prod_{1\le i\le r} (e(\omega_i w)-1)},
\]
where $e(x):=\exp(2\pi\sqrt{-1}x)$, and define a family of polynomials
$P^{(r)}_n(x;\omega_1,\dots,\omega_r)$ in $x$ by the Laurent series expansion
\[
\phi^{(r)}(w;x;\omega_1,\dots,\omega_r)
=
\sum_{n\ge -r} w^n P^{(r)}_n(x;\omega_1,\dots,\omega_r)
\]
valid in a punctured neighborhood of $w=0$; we also set
\[
Q^{(r)}_n(x;\omega_1,\dots,\omega_r)
:=
P^{(r+1)}_n(x;\omega_1,\dots,\omega_r,-1).
\]
By convention, if $n$ is omitted, then $n=0$.
Note that these are related to
the multiple Bernoulli polynomials of \cite{NarukawaA:2004} by
\[
P^{(r)}_n(x;\omega_1,\dots,\omega_r)
=
\frac{(2\pi\sqrt{-1})^n}{(n+r)!} B_{r,n+r}(x|\omega_1,\dots,\omega_r).
\]
Note the special cases
\begin{align}
P^{(0)}(x;) &= 1\notag\\
P^{(1)}(x;\omega) &= \frac{x-\omega/2}{\omega}\notag\\
P^{(2)}(x;\omega_1,\omega_2)
&=
\frac{12(x-\omega_1/2-\omega_2/2)^2-\omega_1^2-\omega_2^2}{24\omega_1\omega_2}
\notag
\end{align}
and
\begin{align}
Q^{(0)}(x;) &= -(x+1/2)\notag\\
Q^{(1)}(x;\omega)
&=
\frac{-12(x-\omega/2+1/2)^2+\omega^2+1}{24\omega}
\notag\\
Q^{(2)}(x;\omega_1,\omega_2)
&=
\frac{
(x-\omega_1/2-\omega_2/2+1/2)
(-4(x-\omega_1/2-\omega_2/2+1/2)^2+\omega_1^2+\omega_2^2+1)
}{24\omega_1\omega_2}.\notag
\end{align}
We will also need a third polynomial
\begin{align}
R^{(r)}(x;\omega_1,\dots,\omega_r)
&=
Q^{(r)}(x;\omega_1,\dots,\omega_r)
+\frac{1}{2}
P^{(r)}(x;\omega_1,\dots,\omega_r)\notag\\
&=
\frac{
Q^{(r)}(x;\omega_1,\dots,\omega_r)
+
Q^{(r)}(x-1;\omega_1,\dots,\omega_r)
}{2}\notag
\end{align}
and note
\begin{align}
R^{(0)}(x;) &= -x\notag\\
R^{(1)}(x;\omega) &= \frac{-12(x-\omega/2)^2+\omega^2-2}{24\omega}\notag\\
R^{(2)}(x;\omega_1,\omega_2) &= 
\frac{
(x-\omega_1/2-\omega_2/2)
(-4(x-\omega_1/2-\omega_2/2)^2+\omega_1^2+\omega_2^2-2)
}{24\omega_1\omega_2}.
\notag
\end{align}

We can then define the {\em hyperbolic gamma function} $\Gammah^{(r)}$ as
follows.  For all integers $r\ge 1$, and $0<\Im(x)<\sum_{1\le i\le
  r}\Im(\omega_i)$,
\[
\Gammah^{(r)}(x;\omega_1,\omega_2,\dots,\omega_r)
:=
\exp(
\PVint_\R\phi^{(r)}(w;x;\omega_1,\dots,\omega_r) \frac{dw}{w}
),
\]
where by the principal value integral $\PVint_\R$ we mean the average of
the integral over two contours agreeing with $\R$ away from 0; one contour
passes to the left of $0$, and the other passes to the right of $0$.  This
differs slightly from the corresponding definitions of hyperbolic gamma
functions in the literature; in particular, we have
\[
\Gammah^{(r)}(x;\omega_1,\dots,\omega_r)
=
S_r(x|\omega_1,\dots,\omega_r)^{(-1)^{r-1}},
\]
where $S_r$ is the multiple sine function (see, for instance,
\cite{NarukawaA:2004}), and
\[
\Gammah^{(2)}(x;\omega_1,\omega_2)
=
G(-\omega_1\sqrt{-1},-\omega_2\sqrt{-1};x-\omega_1/2-\omega_2/2)
\]
in terms of Ruijsenaars' hyperbolic gamma function \cite{RuijsenaarsSNM:1999}.
Since
\[
\PVint_\R
\phi^{(r)}(w;x;\omega_1,\dots,\omega_r)
\frac{dw}{w}
=
\frac{1}{2}
\PVint_\R
\left(\phi^{(r)}(w;x;\omega_1,\dots,\omega_r)
-
\phi^{(r)}(-w;x;\omega_1,\dots,\omega_r)\right)
\frac{dw}{w}
\]
and
\[
\PVint_\R w^k dw = 0
\]
for integers $k<-1$, we can express this as the ordinary integral over $\R$
of
\[
\frac{e(xw)+(-1)^{r-1} e((\sum_{1\le i\le r}\omega_i-x)w)}
     {2w\prod_{1\le i\le r} (e(\omega_iw)-1)}
-
\sum_{1\le k\le (r+1)/2} P^{(r)}_{1-2k}(x;\omega_1,\dots,\omega_r) w^{-2k},
\]
or by symmetry as twice the integral over $[0,\infty)$, thus recovering the
definitions of \cite{RuijsenaarsSNM:1999} and \cite{NarukawaA:2004}.

When $r=1$, we have
\[
\PVint_\R \frac{e(wx)}{(e(w\omega)-1)} \frac{dw}{w}
=
\log(2\sin(\pi x/\omega))
\]
for $0<\Im(x)<\Im(\omega)$ (the branch with value $\log(2)$ at
$x=\omega/2$), which then gives an analytic continuation of
$\Gammah^{(1)}$ to all $x$, $\omega$ such that $\omega\ne 0$:
\[
\Gammah^{(1)}(x;\omega) = 2\sin(\pi x/\omega).
\]
(This is a fairly standard definite integral; it can be shown, for instance
, by moving the contour infinitely far to the left (for $\Im(-x/\omega)>0$;
if $\Im(-x/\omega)<0$, the contour should be moved to the right, and the
case $-x/\omega$ real follows by analytic continuation), and observing that
the resulting sum of residues (compare Theorem \ref{thm:residue_sum} below)
is
\[
\pi\sqrt{-1}(x/\omega)
-\pi\sqrt{-1}/2
-
\sum_{k\ge 1} e(-kx/\omega)/k
=
\log(2\sin(\pi x/\omega)),
\]
as required.)  In general, we have
\[
\Gammah^{(r)}(x+\omega_r;\omega_1,\dots,\omega_r)
=
\Gammah^{(r)}(x;\omega_1,\dots,\omega_r)
\Gammah^{(r-1)}(x;\omega_1,\dots,\omega_{r-1})
\]
in the domain of definition, and thus by induction have a meromorphic
continuation to all $x$.  Since $\Gammah^{(r)}$ was defined via its
logarithm, we in fact have a nearly canonical choice of branch for
$\log\Gammah^{(r)}$; more precisely, for each $\omega_0$ with
$\Im(\omega_0)>0$, there is a unique analytic continuation of
\[
\log\Gammah^{(r)}(x;\omega_1,\dots,\omega_r)
:=
\PVint_\R\phi^{(r)}(w;x;\omega_1,\dots,\omega_r) \frac{dw}{w}
\]
to the set
${\cal C}^{(r)}(\omega_0;\omega_1,\dots,\omega_r)$
obtained from $\C$ by removing the countable unions of rays
\[
\omega_0\R_{\ge 0}+\sum_{1\le i\le r} \omega_i\Z_{\ge 1}
\quad\text{and}\quad
\omega_0\R_{\le 0}+\sum_{1\le i\le r} \omega_i\Z_{\le 0}.
\]
That is, for each zero and pole of $\Gammah^{(r)}$, we cut along a ray in
the direction $\pm\omega_0$, as appropriate.  Also of importance is the
analogous analytic continuation of
\[
\log((x/\omega_0)^{(-1)^r}\Gammah^{(r)}(x;\omega_1,\dots,\omega_r))
\]
to the domain ${\cal C}^{\prime(r)}(\omega_0;\omega_1,\cdots,\omega_r)$
which differs from ${\cal C}^{(r)}(\omega_0;\omega_1,\dots,\omega_r)$ only
in that the ray $\omega_0\R_{\le 0}$ has not been cut; this
continuation exists since the zero/pole at 0 has been cancelled.  We also
by convention take $\Gammah^{(0)}(x;)=-1$, to make the functional equation
valid for $r=1$ as well.  For the functional equation to hold for the
logarithm, we must take
\[
\log\Gammah^{(0)}(x;) := -\pi\sqrt{-1}\sgn(\Im(x/\omega_0)),
\]
defined on the set ${\cal C}^{(0)}(\omega_0;)=\C\setminus \omega_0\R$.

\begin{rem}
Using the fact that
\[
\Gammah^{(r)}(cx;c\omega_1,\dots,c\omega_r)
=
\Gammah^{(r)}(x;\omega_1,\dots,\omega_r)
\]
we can further extend $\Gammah^{(r)}$ to arbitrary $\omega_i$, so long as
there exists a constant $c$ with $\Im(c\omega_i)>0$ for all $i$.  We will
not be using this extension in the sequel, although the ability to rescale
within the upper half-plane will be quite useful in the proofs.
\end{rem}

In \cite[Prop. 5]{NarukawaA:2004}, Narukawa derived a product expansion for
$\Gammah^{(r)}$, based on the observation that the integral over a
suitably chosen sequence of contours $\Im(w)=a$ can be made to tend to 0 as
$a\to\infty$, and thus the integral expands as a sum of residues.  We will
need a more precise form of the bound on the integral.

\begin{lem}\label{lem:hyperb_asympt}
Fix $\epsilon>0$, and let $a$, $\Im(x)$, $\omega_i$ range over the domain
$0<\Im(x)<\sum_i\Im(\omega_i)$, and for $1\le i\le r$, $\Im(\omega_i)>0$
and $|a-\Im(-n/\omega_i)|>\epsilon$, all integers $n\ge 0$.
Then
\[
\left|\int_{\Im(w)=a} \phi^{(r)}(w;x;\omega_1,\dots,\omega_r) \frac{dw}{w}\right|
\le
a^{-1}
\left[
\prod_i C(\epsilon|\omega_i|)^{-1}
\right]
\left[
\frac{\exp(-2\pi a\Re(x))}
     {2\pi\Im(x)}
+
\frac{\exp(-2\pi a\Re(x-\sum_i\omega_i))}
     {2\pi\Im(\sum_i\omega_i-x)}
\right],
\]
where $C(x) = \min_{d(y,\Z)\ge x} |e(y)-1|.$
In particular, the integral is uniformly $O(\exp(-2\pi a\Re(x)))$ over any
compact subset of the domain.
\end{lem}

\begin{proof}
We have the estimates
\begin{align}
|e(wx)|  &= \exp(-2\pi \Re(w)\Im(x))\exp(-2\pi a\Re(x))\notag\\
|w|^{-1} &= 1/\sqrt{\Re(w)^2+a^2}\le 1/a\notag\\
|e(w\omega_i)-1|^{-1} &\le C(\epsilon|\omega_i|)^{-1}\notag\\
|e(w\omega_i)-1|^{-1} &\le C(\epsilon|\omega_i|)^{-1}
\exp(2\pi\Re(w)\Im(\omega_i))\exp(2\pi a \Re(\omega_i))\notag
\end{align}
where we note that $|w\omega_i-n|>\epsilon|\omega_i|$.
If we use the third estimate for $\Re(w)>0$ and the fourth
estimate for $\Re(w)<0$, we obtain the stated bound.
\end{proof}

Narukawa then uses the fact that for $\Re(x)>\max(0,\Re(\sum_i\omega_i))$,
the above bound tends to 0 as $a\to \infty$; for our purposes, it is more
convenient to fix $a$ and obtain an asymptotic series.

\begin{thm}\label{thm:residue_sum}
Let $a$, $\arg(x)$, $\omega_i$ range over the domain $a>0$,
$\Im(e(-\arg(x)/2\pi)\omega_i)>0$ for $1\le i\le r$.  Then as $x\to\infty$,
\begin{align}
-\pi\sqrt{-1}P^{(r)}(x;\omega_1,\dots,\omega_r)
-
\sum_{0<\Im(e(\arg(x)/2\pi)y)\le a}
 \Res_{w=y}(\phi^{(r)}(w;x;\omega_1,\dots,\omega_r) \frac{dw}{w})
+\log&\Gammah^{(r)}(x;\omega_1,\dots,\omega_r)
\notag\\
&
{}=O(\exp(-2\pi a |x|)),\notag
\end{align}
uniformly over any compact subset of the domain.
\end{thm}

\begin{proof}
Clearly replacing $\arg(x)$ by $\arg(x-b\sum_i\omega_i)$ for $0<b<1$ will
have no effect on the validity of the bound.  We may thus take
\[
x=e(\tau)|y|+b\sum_i\omega_i.
\]
Since
\[
\log\Gammah^{(r)}(e(\phi)|y| + b\sum_i\omega_i;\omega_1,\dots,\omega_r)
=
\log\Gammah^{(r)}(|y| + b\sum_i
e(-\phi)\omega_i;e(-\phi)\omega_1,\dots,e(-\phi)\omega_r),
\]
we find that uniformity in $a,\tau,\omega_i$ will follow from uniformity in
$a,\omega_i$ with $\tau=0$.

Thus assume $\Im(x)=b\sum_i\Im(\omega_i)$ for $0<b<1$.  For every point in
the domain, there exists $a'\ge a$ and $\epsilon>0$ such that the previous
lemma applies, and such that there are no poles with imaginary part in
$(a,a']$; by compactness, we can cover the domain by a finite number of
such choices.  Since $O(\exp(-2\pi a'\Re(x))) = O(\exp(-2\pi a\Re(x)))$, it
suffices to consider the case $a'=a$.  The result follows by residue
calculus.
\end{proof}

\begin{rem}
Note that the cut lines for $\log\Gammah^{(r)}$ can be taken along any
direction $\omega_0$ in the convex cone generated by
$\omega_1$,\dots,$\omega_r$, and the argument is valid for $x$ in the
complement of the cones
\[
\omega_1+\cdots+\omega_r + \R_{\le 0}\langle \omega_1,\dots,\omega_r\rangle
\quad\text{and}\quad
-\R_{\le 0}\langle \omega_1,\dots,\omega_r\rangle,
\]
and not approaching $\infty$ parallel to the boundary of the cone.
\end{rem}

In particular, we obtain the following estimate.

\begin{cor}\label{cor:hyperb_asympt}
Let $a$, $\arg(x)$, $\omega_i$ range over the domain
$0<a<\min_i\Im(-e(\arg(x)/2\pi)/\omega_i)$.  Then as $x\to\infty$, we have the
estimates
\begin{align}
-\pi\sqrt{-1}P^{(r)}(x;\omega_1,\dots,\omega_r)
+\log\Gammah^{(r)}(x;\omega_1,\dots,\omega_r)
&=
O(\exp(-2\pi a |x|)),\notag\\
\pi\sqrt{-1}P^{(r)}(-x;\omega_1,\dots,\omega_r)
+\log\Gammah^{(r)}(-x;\omega_1,\dots,\omega_r)
&=
O(\exp(-2\pi a |x|)),\notag
\end{align}
uniformly on compact subsets.
\end{cor}

\begin{proof}
The first estimate follows immediately from the theorem.
The second estimate then follows using the facts
\[
P^{(r)}(x;\omega_1,\dots,\omega_r)
=
(-1)^r P^{(r)}(\sum_i\omega_i-x;\omega_1,\dots,\omega_r)
\]
and
\[
\log\Gammah^{(r)}(x;\omega_1,\dots,\omega_r)
=
(-1)^{r-1}\log\Gammah^{(r)}(\sum_i\omega_i-x;\omega_1,\dots,\omega_r)
\]
\end{proof}

\begin{rem}
For $r=2$, this result is essentially due to Ruijsenaars
\cite[App. 2]{RuijsenaarsSNM:1999}.  In general, on the domain
$\Im(-e(\arg(x)/2\pi)/\omega_i)>0$, it should be possible to improve the error
term to
\[
O(|x|^{r-1}\exp(-2\pi |x|\min_i\Im(-e(\arg(x)/2\pi)/\omega_i)),
\]
by taking $a=\min_i \Im(-e(\arg(x)/2\pi)/\omega_i)$ and bounding the leading
residues.  In particular, if the poles with minimum imaginary part are all
simple in the given compact subset of parameter space, then one easily has
the bound
\[
O(\exp(-2\pi |x|\min_i\Im(-e(\arg(x)/2\pi)/\omega_i)))
\]
on their residues, and thus on the error term.
\end{rem}

We will also need an extension of this to the case that one of the moduli
tends to 0.  For convenience, write
\[
\gamma^{(r)}_h(x;\omega_1,\dots,\omega_r)
=
-\pi\sqrt{-1}P^{(r)}(x;\omega_1,\dots,\omega_r)
+\log\Gammah^{(r)}(x;\omega_1,\dots,\omega_r).
\]

\begin{thm}\label{thm:gasympt:xbig_vsmall}
Let $a$, $\arg(x)$, $\omega_i$,$\alpha$,$\beta$ range over the domain
\[
0<a<\min_{1\le i\le r-1}\Im(-e(\arg(x)/2\pi)/\omega_i),\quad
\Im(e(-\arg(x)/2\pi)/\omega_r)>0.
\]
Then as $x\to\infty$, $v\to 0^+$,
\begin{align}
 \gamma_h^{(r)}(x+v\alpha\omega_r;\omega_1,\dots,\omega_{r-1},v\omega_r)
-\gamma_h^{(r)}(x+v\beta\omega_r;\omega_1,\dots,\omega_{r-1},v\omega_r)
-(\alpha-\beta)\gamma^{(r-1)}_h&(x;\omega_1,\dots,\omega_{r-1})\notag\\
&{}=
O(v\exp(-2\pi a|x|)),\notag
\end{align}
uniformly over compact subsets of the domain.
\end{thm}

\begin{proof}
For $r=1$, this is immediate, so take $r\ge 2$.  As before, we may restrict
our attention to the case $0<\Im(x)<\Im(\sum_{1\le i\le r-1}\omega_i)$.
Then the left-hand side can be expressed as
\[
\int_{\Re(w)=a}
\left(
\frac{e(v\alpha\omega_r w)-e(v\beta\omega_r w)}{e(v\omega_r w)-1}
-
(\alpha-\beta)
\right)
\phi^{(r-1)}(w;x;\omega_1,\dots,\omega_{r-1})
\frac{dw}{w}
\]
The desired estimate then follows as in Lemma \ref{lem:hyperb_asympt},
using the fact that as $v\to 0^+$,
\[
\frac{1}{v}
\left(
\frac{e(v\alpha\omega_r w)-e(v\beta\omega_r w)}{e(v\omega_r w)-1}
-
(\alpha-\beta)
\right)
\]
is $O(\exp(\epsilon |\Re(w)|))$ for any $\epsilon>0$, uniformly in $w$ and on
compact subsets of parameter space.
\end{proof}

\begin{cor}\label{cor:hyper_to_rat}
Let $x$, $\omega_1$,\dots,$\omega_r$, $\alpha$, $\beta$ range over the domain
$\Im(\omega_i)>0$, $x\in {\cal
  C}^{(r-1)}(\omega_r;\omega_1,\dots,\omega_{r-1})$.
Then as $v\to 0^+$, we have the estimate
\begin{align}
 \log
\frac{\Gammah^{(r)}(x+v\alpha\omega_r;\omega_1,\dots,\omega_{r-1},v\omega_r)}
{\Gammah^{(r)}(x+v\beta\omega_r;\omega_1,\dots,\omega_{r-1},v\omega_r)}
-(\alpha-\beta)\log\Gammah^{(r-1)}(x;\omega_1,\dots,\omega_{r-1})
=
O(v),\notag
\end{align}
uniformly over compact subsets of the domain.
\end{cor}

\begin{proof}
Since
\[
{\cal
  C}^{(r-1)}(\omega_r;\omega_1,\dots,\omega_{r-1})
=
{\cal C}^{(r)}(\omega_r;\omega_1,\dots,\omega_r),
\]
there is a canonical choice of branch for the left-hand side, and the proof
of Theorem \ref{thm:gasympt:xbig_vsmall} gives the desired estimate on
compact subsets of the complement of the cones
\[
\omega_1+\cdots+\omega_r + \R_{\le 0}\langle \omega_1,\dots,\omega_r\rangle
\quad\text{and}\quad
-\R_{\le 0}\langle \omega_1,\dots,\omega_r\rangle.
\]
Since the estimate is consistent with the functional equation for
$\Gammah^{(r)}$, the result follows.
\end{proof}

\begin{rem}
In both cases, this can be extended to a uniform asymptotic series
\[
 \log
\frac{\Gammah^{(r)}(x+v\alpha\omega_r;\omega_1,\dots,\omega_{r-1},v\omega_r)}
{\Gammah^{(r)}(x+v\beta\omega_r;\omega_1,\dots,\omega_{r-1},v\omega_r)}
=
\sum_{1\le k\le n}
\frac{B_k(\alpha)-B_k(\beta)}{k!}
\left(\omega_r v\frac{d}{dx}\right)^{k-1}
\log\Gammah^{(r-1)}(x;\omega_1,\dots,\omega_{r-1})
+
O(v^n),
\]
with error $O(v^n \exp(-2\pi a|x|))$ if also $x\to\infty$; here $B_k(x)$ is
the $k$th (ordinary) Bernoulli polynomial.
\end{rem}

In fact, with care, we can give an estimate valid on the larger domain
${\cal C}^{\prime(r-1)}(\omega_r;\omega_1,\dots,\omega_{r-1})$, and thus in
particular in a neighborhood of 0.  The point is that we can identify the
zeros and poles of $\Gammah^{(r)}$ and $\Gammah^{(r-1)}$ that give rise
to the cut line $\omega_r\R_{\le 0}$ above, and using the ordinary gamma
function, cancel them out.  The resulting asymptotics can then still be
computed via Stirling's formula.  It will be convenient to use a slightly
renormalized form of the ordinary gamma function; we define
\[
\Gammar(x;\omega) := \frac{\Gamma(x/\omega)}{\sqrt{2\pi}},
\]
with the convention that $\omega=1$ if omitted, and the usual convention on
multiple arguments.  As a justification for this convention, note that the
reflection identity for the ordinary gamma function becomes
\[
\Gammar(x,\omega-x;\omega) = \Gammah^{(1)}(x;\omega)^{-1}.
\]
We also let $\log\Gammar(x;\omega)$ denote the standard branch on
$\C\setminus \omega\R_{\le 0}$, and note the following version of
Stirling's formula.

\begin{thm}
Let $x\notin \R_{\le 0}$. Then for all $m\ge 1$,
\[
\log\Gammar(x)
=
(x\log(x)-x) - \frac{1}{2}\log(x)
+
\sum_{1\le i<m} \frac{B_{2i} x^{1-2i}}{2i(2i-1)}
+
O(d(x,\R_{\le 0})^{1-2m}),
\]
uniformly in $x$ as $d(x,\R_{\le 0})\to\infty$.  More generally, for
$\alpha$ ranging over any compact subset of $\C$,
\[
\log\Gammar(x+\alpha)
=
(x\log(x)-x) +B_1(\alpha)\log(x)
+
\sum_{1\le i<m} \frac{(-1)^{i+1} B_{i+1}(\alpha)}{i(i+1) x^i}
+
O(d(x,\R_{\le 0})^{-m}),
\]
uniformly as $d(x,\R_{\le 0})\to\infty$.
\end{thm}

\begin{proof}
The claim for general $\alpha$ follows from the claim for $\alpha=0$ by
straightforward algebraic manipulation.  For $\alpha=0$, we observe that by
\cite[\S 8.4]{OlverFWJ:1974} or
\cite[Thm. 1.4.2]{AndrewsGE/AskeyR/RoyR:1999}, the error term is
\[
\frac{B_{2m} x^{1-2m}}{2m(2m-1)}
-
\frac{1}{2m}
\int_0^\infty
\frac{B_{2m}(\{t\})}{(x+t)^{2m}}
dt.
\]
The first term certainly has the correct asymptotics; for the second term,
we have
\[
\frac{1}{2m}
\int_0^\infty
\frac{B_{2m}(\{t\})}{(x+t)^{2m}}
dt
=
O(\int_0^\infty |x+t|^{-2m} dt).
\]
For $\Re(x)\ge 0$, $|x+t|^{-2m}\le (|x|^2+t^2)^{-m}$, and thus the
integral has order $|x|^{1-2m}$.  For $\Re(x)\le 0$, the integral is still
bounded above by
\[
\int_{-\infty}^{\infty} (\Im(x)^2+t^2)^{-m} dx = O(|\Im(x)|^{1-2m}).
\]
\end{proof}

\begin{rem}
More generally, \cite[Ex. 8.4.4]{OlverFWJ:1974} gives the error term
\[
\frac{(-1)^{m+1} B_{m+1}(\alpha) x^{-m}}{2m(2m+1)}
-
\frac{1}{m+1}
\int_0^\infty
\frac{B_{m+1}(\{t-\alpha\})}{(x+t)^{m+1}}
dt
\]
for $0\le \alpha\le 1$.
\end{rem}

\begin{cor}
Let $x,\omega$,$\alpha$,$\beta$ range over the domain $\omega\ne 0$, $x\in
\C\setminus \omega\R_{\le 0}$.  Then as $v\to 0^+$,
\begin{align}
\log
\frac{\Gammar(x+v\alpha\omega;v\omega)}
     {\Gammar(x+v\beta\omega;v\omega)}
={}&
(B_1(\alpha)-B_1(\beta))\log(x/v\omega)
+
\sum_{2\le i<m} \frac{(-1)^i (B_i(\alpha)-B_i(\beta)) (x/v\omega)^{1-i}}{i(i-1)}\notag\\
&+
O(d(x,\omega\R_{\le 0})^{-m}v^m),\notag
\end{align}
uniformly in $x$ and over compacta in $\omega$, $\alpha$, $\beta$.
\end{cor}

\begin{thm}\label{thm:hyper_to_rat}
Let $x$, $\omega_1$,\dots,$\omega_r$, $\alpha$,$\beta$ range over the domain
$\Im(\omega_i)>0$, $x\in {\cal
  C}^{\prime (r-1)}(\omega_r;\omega_1,\dots,\omega_{r-1})$.
Then as $v\to 0^+$, we have the estimate
\begin{align}
 \log
\frac{\Gammah^{(r)}(x+v\alpha\omega_r;\omega_1,\dots,\omega_{r-1},v\omega_r)}
     {\Gammah^{(r)}(x+v\beta\omega_r;\omega_1,\dots,\omega_{r-1},v\omega_r)}
&{}-
(-1)^r
\log\frac{\Gammar(x+v\alpha\omega_r;v\omega_r)}
         {\Gammar(x+v\beta\omega_r;v\omega_r)}\notag\\
&{}-(\alpha-\beta)
\left(
\log\Gammah^{(r-1)}(x;\omega_1,\dots,\omega_{r-1})
-(-1)^r\log(x/v\omega_r)
\right)
=
O(v),\notag
\end{align}
uniformly over compact subsets of the domain.
\end{thm}

\begin{proof}
We first observe that
\[
\log\Gammah^{(r)}(x+v\alpha\omega_r;\omega_1,\dots,\omega_{r-1})
-(-1)^r\log\Gammar(x+v\alpha\omega_r;v\omega_r)
\]
and
\[
\log\Gammah^{(r-1)}(x;\omega_1,\dots,\omega_{r-1})
-(-1)^r\log(x/v\omega_r)
\]
are analytic on the given domain, and thus the overall left-hand side is
analytic.  Moreover, the stated estimate holds on the smaller domain ${\cal
  C}^{(r-1)}(\omega_r;\omega_1,\dots,\omega_{r-1})$.  Using the functional
equation, we may immediately extend this to the full domain ${\cal
  C}^{\prime(r-1)}(\omega_r;\omega_1,\dots,\omega_{r-1})$, except when
$r=2$, where the point $x=0$ must still be excluded (since the only points
related to $x=0$ via the functional equation are also in the cut set).  But
in that case, we may simply use Cauchy's theorem to deduce a uniform
estimate on a neighborhood of $0$ from a uniform estimate on the boundary
of the neighborhood.
\end{proof}

\begin{rem}
More generally one has an asymptotic series in which the coefficient of the
$k$th term depends on the $(k-1)$-st derivative of
\[
\log\Gammah^{(r-1)}(x;\omega_1,\dots,\omega_{r-1})
-(-1)^r\log(x/v\omega_r).
\]
\end{rem}

If $r=2$, we have
\[
\lim_{x\to 0} \log\Gammah^{(2)}(x;\omega_1,\omega_2)-\log\Gammar(x;\omega_2)
=
\frac{\log(\omega_1)-\log(\omega_2)-\log(2\pi)}{2},
\]
which gives rise to a nice corollary by taking $x=\beta=0$ above.

\begin{cor}\cite{RuijsenaarsSNM:1999}\label{cor:hyper_to_rat2}
As $v\to 0^+$,
\[
\log\Gammah^{(2)}(v\alpha\omega_2;\omega_1,v\omega_2)
-\log\Gammar(\alpha)
-(\alpha-1/2)\log(2\pi v\omega_2/\omega_1)
=
O(v)
\]
uniformly over compact subsets of the region
$\Im(\omega_1),\Im(\omega_2)>0$, $\alpha\in \C$.
\end{cor}

Now, consider the elliptic gamma function, defined as
\[
\Gammae^{(r)}(z;p_1,p_2,\dots,p_r)
=
\prod_{0\le k_1,k_2,\dots,k_r}
(1-p_1^{k_1+1}p_2^{k_2+1}\cdots p_r^{k_r+1}/z)
(1-p_1^{k_1}p_2^{k_2}\cdots p_r^{k_r} z)^{(-1)^{r-1}}.
\]
For $|p_1p_2\cdots p_r|<|z|<1$, we have
\[
\log\Gammae^{(r)}(z;p_1,p_2,\dots,p_r)
=
\sum_{1\le k}
\frac{(-1)^r z^k-(p_1p_2\cdots p_r/z)^k}
     {k\prod_{1\le i\le r} (1-p_i^k)};
\]
this, then, for any $\omega_0$ with $\Im(\omega_0)>0$, defines a branch of
$\log\Gammae^{(r)}$ on the region ${\cal
  C}^{(r)}_e(\omega_0;p_1,\dots,p_r)$ obtained from $\C^*=\C\setminus
\{0\}$ by removing the countable union of logarithmic spirals
\[
p_1^{\Z_{\ge 1}}p_2^{\Z_{\ge 1}}\cdots p_r^{\Z_{\ge 1}}e(\omega_0\R_{\ge
  0})
\quad\text{and}\quad
p_1^{\Z_{\le 0}}p_2^{\Z_{\le 0}}\cdots p_r^{\Z_{\le 0}}e(\omega_0\R_{\le
  0}).
\]
(We define a region ${\cal C}^{\prime(r)}_e$ analogously.)  And, of course,
we have the functional equation
\[
\log\Gammae^{(r)}(p_rz;p_1,p_2,\dots,p_r)
-
\log\Gammae^{(r)}(z;p_1,p_2,\dots,p_r)
=
\log\Gammae^{(r-1)}(z;p_1,p_2,\dots,p_{r-1}),
\]
with
\[
\log\Gammae^{(0)}(z;)=\log(-1/z),\quad \log\Gammae^{(0)}(-1;)=0.
\]

Narukawa \cite[Theorem 14]{NarukawaA:2004} gives the following ``product''
expansion of the
elliptic gamma function in terms of the hyperbolic gamma function:
\begin{align}
-2\pi\sqrt{-1}Q^{(r)}(x;\omega_1,\dots,\omega_r)
+
{}&\log\Gammae^{(r)}(e(x);e(\omega_1),\dots,e(\omega_r))\notag\\
&{}=
\sum_{0\le k}
-\pi\sqrt{-1}P^{(r)}(x+k+1;\omega_1,\dots,\omega_r)
+\log\Gammah^{(r)}(x+k+1;\omega_1,\dots,\omega_r)\notag\\
&+\sum_{0\le k}
\pi\sqrt{-1}P^{(r)}(x-k;\omega_1,\dots,\omega_r)
+\log\Gammah^{(r)}(x-k;\omega_1,\dots,\omega_r)\notag
\end{align}
Note that each term in the infinite sums converges uniformly exponentially
to 0 as $k\to\infty$, so the sums converge uniformly and absolutely.

Using this expansion and the asymptotics of the hyperbolic gamma
function, we obtain the following estimates for the elliptic gamma function.
First, the hyperbolic limit $p_1,\dots,p_r\to 1$.

\begin{prop}\label{prop:ell_to_hyper2}
Let $A$, $\omega_1,\dots,\omega_r$, $x$ range over the domain
\[
0<A<\min_{1\le i\le r} \Im(-1/\omega_i)-|\Im(x/\omega_i)|.
\]
Then as $v\to 0^+$, we have the estimate
\[
-2\pi\sqrt{-1}R^{(r)}(x;v\omega_1,\dots,v\omega_r)
+\log\Gammae^{(r)}(e(x);e(v\omega_1),\dots,e(v\omega_r))
-\log\Gammah^{(r)}(x;v\omega_1,\dots,v\omega_r)
=
O(\exp(-2\pi A/v)),
\]
uniform over compact subsets of the domain.
\end{prop}

\begin{rem}
If we do not care about the choice of $A$, the
constraint on the domain is simply
\[
|\Im(x/\omega_i)|<\Im(-1/\omega_i), 1\le i\le r;
\]
this is a parallelogram, two of the vertices of which are $\pm 1$.
\end{rem}

\begin{prop}\label{prop:ell_to_hyper1}
Let $A$, $\omega_1,\dots,\omega_r$, $x$ range over the domain
\[
0<A<\min_{1\le i\le r} \min(\Im(-x/\omega_i),\Im((x-1)/\omega_i)).
\]
Then as $v\to 0^+$, we have the estimate
\[
-2\pi\sqrt{-1}Q^{(r)}(x;v\omega_1,\dots,v\omega_r)
+\log\Gammae^{(r)}(e(x);e(v\omega_1),\dots,e(v\omega_r))
=
O(\exp(-2\pi A/v)),
\]
uniform over compact subsets of the domain.
\end{prop}

\begin{proof}
Using the ``product'' expansion and the identity
\[
\log\Gammah^{(r)}(z;v\omega_1,\dots,v\omega_r)
=
\log\Gammah^{(r)}(z/v;\omega_1,\dots,\omega_r),
\]
we can express each left-hand side as a sum over functions to which
Corollary \ref{cor:hyperb_asympt} applies, giving the desired uniform
asymptotics.
\end{proof}

Using this limit, we can obtain the following bound.

\begin{prop}\label{prop:ell_bound_hyper}
Fix a compact subset $S$ of the set of $\omega_i$ such that
$0<\Im(\omega_i)$, and constants $0<\epsilon,C_1,C_2$.  Then as
$v\to 0^+$, we have the following estimate, uniform over the product of $S$
with the region $-1-vC_1\le \Re(x)\le vC_1$, $|\Im(x)|\le v C_2$, apart
from a hole of radius $\epsilon v$ around every pole of the left-hand side:
\[
\Gammae^{(r)}(e(x);e(v\omega_1),\dots,e(v\omega_r))^{\pm 1}
=
O(e(\pm Q^{(r)}(x;v\omega_1,\dots,v\omega_r))).
\]
\end{prop}

\begin{proof}
We consider the $+$ case; the $-$ case is completely analogous.  Choose
$1/2<D<1$ and $0<a<\min_i\Im(-1/\omega_i)$.  Proposition
\ref{prop:ell_to_hyper2} gives
\[
e(-Q^{(r)}(x;v\omega_1,\dots,v\omega_r))
\Gammae^{(r)}(e(x);e(v\omega_1),\dots,e(v\omega_r))
=
O(
e(\frac{1}{2}P^{(r)}(x;\omega_1,\dots,\omega_r))
\Gammah^{(r)}(x;v\omega_1,\dots,v\omega_r)
)
\]
away from the poles, in the subregion $-D\le \Re(x)\le vC_1$.  Since
\[
e(\frac{1}{2}P^{(r)}(x;v\omega_1,\dots,v\omega_r))
\Gammah^{(r)}(x;v\omega_1,\dots,v\omega_r)
\]
is uniformly bounded in that region (we have excluded neighborhoods of the
poles, and it converges uniformly to 1 for $-\Re(x)/v$ large), we have the
uniform estimate
\[
e(-Q^{(r)}(x;v\omega_1,\dots,v\omega_r))
\Gammae^{(r)}(e(x);e(v\omega_1),\dots,e(v\omega_r))
=
O(1)
\]
in this region.  A similar application of Proposition \ref{prop:ell_to_hyper2}
gives
\begin{align}
e(-Q^{(r)}(x+1;v\omega_1,\dots,v\omega_r))
&\Gammae^{(r)}(e(x);e(v\omega_1),\dots,e(v\omega_r))
\notag\\
&=
O(
e(\frac{1}{2}P^{(r)}(x+1;v\omega_1,\dots,v\omega_r))
\Gammah^{(r)}(x+1;v\omega_1,\dots,v\omega_r)
)
\notag
\end{align}
uniformly over the subset $-vC_1\le \Re(x)+1\le D$.  Since
\[
Q^{(r)}(x+1;v\omega_1,\dots,v\omega_r)
=
Q^{(r)}(x;v\omega_1,\dots,v\omega_r)
-
P^{(r)}(x+1;v\omega_1,\dots,v\omega_r),
\]
we find
\begin{align}
e(-Q^{(r)}(x;v\omega_1,\dots,v\omega_r))
&\Gammae^{(r)}(e(x);e(v\omega_1),\dots,e(v\omega_r))
\notag\\
&=
O(
e(-\frac{1}{2}P^{(r)}(x+1;v\omega_1,\dots,v\omega_r))
\Gammah^{(r)}(x+1;v\omega_1,\dots,v\omega_r)
)
=
O(1)
\notag
\end{align}
on this region as well.
\end{proof}

Similarly, we can obtain asymptotics of $\Gammae^{(r)}$ in the
``rational'' limit $p_r\to 1$, $p_1$,\dots $p_{r-1}$ fixed.

\begin{thm}\label{thm:ell_to_rat1}
Let $x$, $p_1$,\dots,$p_{r-1}$, $\omega_r$, $\alpha$, $\beta$ range over
the domain $0<|p_1|,\dots,|p_{r-1}|<1$, $\Im(\omega_r)>0$, and
\[
x\in e^{-1}({\cal C}^{(r-1)}_e(\omega_r;p_1,\dots,p_{r-1}))
\cup
(\omega_r\R_{\ge 0}\cap e^{-1}({\cal C}^{\prime
  (r-1)}_e(\omega_r;p_1,\dots,p_{r-1})).
\]
Then as $v\to 0^+$, we have the estimate
\begin{align}
 \log
\frac{\Gammae^{(r)}(e(x+v\alpha\omega_r);p_1,\dots,p_{r-1},e(v\omega_r))}
     {\Gammae^{(r)}(e(x+v\beta\omega_r);p_1,\dots,p_{r-1},e(v\omega_r))}
&{}-
(-1)^r
\log\frac{\Gammar(x+v\alpha\omega_r;v\omega_r)}
         {\Gammar(x+v\beta\omega_r;v\omega_r)}\notag\\
{}-(\alpha-\beta)&
\left(
\log\Gammae^{(r-1)}(e(x);p_1,\dots,p_{r-1})
-(-1)^r\log(x/v\omega_r)
\right)
=
O(v),\notag
\end{align}
uniformly over compact subsets of the domain.
\end{thm}

\begin{proof}
If we expand the elliptic gamma functions via the product representation
and group corresponding terms, we find by Corollary \ref{cor:hyper_to_rat}
that all but one term of the result is uniformly $O(v)$.  Moreover,
it follows from Theorem \ref{thm:gasympt:xbig_vsmall} that the coefficient
of $v$ in the estimates is exponentially small as $k\to \infty$, and thus
the error terms are summable.  The only surviving term can be estimated
using Theorem \ref{thm:hyper_to_rat}, giving the desired result.
\end{proof}

\begin{cor}\label{cor:ell_to_rat2}
As $v\to 0^+$,
\[
\log
\frac{
\Gammae^{(2)}(e(\alpha v\omega);p,e(v\omega))
}{
\Gammar(\alpha)
}
=
\frac{\pi\sqrt{-1}}{12v\omega}
+
(\alpha-1/2)
\log(2\pi v\omega (p;p)^2/\sqrt{-1})
+
O(v),
\]
uniformly over compact subsets of the region
$\Im(\omega)>0$, $0<|p|<1$, $\alpha\in \C$.
\end{cor}

Similarly, one has the following.

\begin{cor}\label{cor:ell_to_rat3}
As $v\to 0^+$,
\[
\frac{
\Gammae^{(r)}(e(\alpha v\omega) z;p_1,\dots,p_{r-1},e(v\omega))
}{
\Gammae^{(r)}(e(\beta v\omega) z;p_1,\dots,p_{r-1},e(v\omega))
}
=
\Gammae^{(r-1)}(z;p_1,\dots,p_{r-1})^{\alpha-\beta}
(1+O(v))
\]
uniformly over compact subsets of the domain $0<|p_1|,\dots,|p_{r-1}|<1$,
$\Im(\omega)>0$, $|p_1\cdots p_{r-1}|<z\le \min_{1\le i\le r-1}
|p_i|^{-1}$,  $z\ne 1$, $\alpha,\beta\in \C$.
\end{cor}

\begin{cor}\label{cor:ell_bound_rat}
As $v\to 0^+$, the function
\[
\frac{
\Gammae^{(r)}(e(v\alpha\omega) z;p_1,\dots,p_r,e(v\omega))
}{
\Gammae^{(r)}(e(v\beta\omega) z;p_1,\dots,p_r,e(v\omega))
}
\]
is uniformly bounded over compact subsets of the domain
$0<|p_1|,\dots,|p_{r-1}|<1$, $\Im(\omega)>0$, $|p_1\cdots p_{r-1}|<z<1$,
$\alpha,\beta\in\C$.  On compact subsets of the domain
$0<|p_1|,\dots,|p_{r-1}|<1$, $\Im(\omega)>0$,
$|p|<z\le 1$, $\alpha,\beta\in \C$, it is uniformly
\[
O(v^{\min(\Re(\alpha-\beta),0)}).
\]
\end{cor}

There are corresponding estimates for the trigonometric gamma function
\[
\Gammat(z;q) := \prod_{0\le k} (1-q^k z)^{-1} = \frac{1}{(z;q)_\infty},
\]
with the corresponding analytic continuation of its logarithm (the branch
with $\log\Gammat(0;q)=0$).

\begin{lem}\cite{KoornwinderTH:1990}\label{lem:trig_to_rat0}
Let $z$, $\omega$, $\alpha$, $\beta$ range over the domain $\Im(\omega)>0$,
$0\le |z|<1$.  Then as $v\to 0$, we have the estimate
\[
\log
\frac{\Gammat(e(v\alpha\omega)z;e(v\omega))}
     {\Gammat(e(v\beta\omega)z;e(v\omega))}
-
(\alpha-\beta)\log(1-z)
=
O(vz),
\]
uniformly over compact subsets of the domain.
\end{lem}

\begin{proof}
Indeed, for $|z|<1$, the left-hand side is given by the sum
\[
\sum_{k\ge 1}
\frac{z^k}{k}
\left(
\frac{e(kv\alpha\omega)-e(kv\beta\omega)}
     {1-e(kv\omega)}
+
\alpha-\beta
\right)
\]
and the quantity in parentheses is uniformly $O(v\exp(\epsilon k))$ for all
$\epsilon>0$.
\end{proof}

This gives rise to a trigonometric analogue of Theorem
\ref{thm:ell_to_rat1}.

\begin{thm}
For $r>1$, let $z$, $p_1$,\dots,$p_{r-1}$, $\omega$, $\alpha$, $\beta$
range over the domain $0\le |p_1|,\dots,|p_{r-1}|<1$, $\Im(\omega)>0$, and
$z\in {\cal C}^{\prime(r-1)}_e(\omega;p_1,\dots,p_{r-1})$.
Then as $v\to 0^+$, we have the estimate
\begin{align}
\log
\frac{\Gammae^{(r)}(e(v\alpha\omega)z;p_1,\dots,p_{r-1},e(v\omega))}
     {\Gammae^{(r)}(e(v\beta\omega)z;p_1,\dots,p_{r-1},e(v\omega))}
&{}-
(-1)^r
\log\frac{\Gammat(e(v\alpha\omega)z;e(v\omega))}
         {\Gammat(e(v\beta\omega)z;e(v\omega))}\notag\\
{}-(\alpha-\beta)&
\left(
\log\Gammae^{(r-1)}(z;p_1,\dots,p_{r-1})
-(-1)^r\log(1-z)
\right)
=
O(v),\notag
\end{align}
uniformly over compact subsets of the domain.
\end{thm}

\begin{proof}
For $|p_1p_2\cdots p_{r-1}|<z<\min_i |p_i|^{-1}$, this follows immediately from
the expansion
\[
\log\Gammae^{(r)}(z;p_1,\dots,p_{r-1},q)
=
-
\sum_{0\le k_1,k_2,\dots,k_{r-1}}
\Gammat(p_1^{k_1+1}\cdots p_{r-1}^{k_{r-1}+1}/z;q)
+(-1)^{r-1}\Gammat(p_1^{k_1}\cdots p_{r-1}^{k_{r-1}}z;q)
\]
and the asymptotics of $\Gammat$.

The general case follows from the functional equation.
\end{proof}

Comparing this to Theorem \ref{thm:ell_to_rat1} gives the following result,
a uniform version of the results of the appendices of
\cite{KoornwinderTH:1990}.

\begin{cor}\label{cor:trig_to_rat}
Let $x$, $\omega$, $\alpha$, $\beta$ range over the domain $\Im(\omega)>0$,
\[
x\in \bigl(\C\setminus (\Z+\omega\R_{\ge 0})\bigr)\cup \omega \R_{\ge 0}.
\]
Then as $v\to 0^+$, we have the estimate
\[
\log
\frac{\Gammat(e(x+v\alpha\omega);e(v\omega))}
     {\Gammat(e(x+v\beta\omega);e(v\omega))}
-
\log
\frac{\Gammar(x+v\alpha\omega;v\omega)}
     {\Gammar(x+v\beta\omega;v\omega)}
-
(\alpha-\beta)
\bigl(\log(1-e(x))-\log(x/v\omega)\bigr)
=
O(v),
\]
uniformly over compact subsets of the domain.
\end{cor}

\begin{rem}
Note that the validity of the theorem for $r=2$, $|p|<z<|p|^{-1}$ is enough
to give the corollary in general, and in turn give the theorem in general,
without using the functional equation.  This also implies that Theorem
\ref{thm:ell_to_rat1} and its corollaries continue to hold even without the
constraint $0<|p_i|$ on the domain, and further implies that Lemma
\ref{lem:trig_to_rat0} holds on the domain $z\notin e(\omega\R_{\ge 0})$.
\end{rem}

\begin{cor}\cite{KoornwinderTH:1990}
Let $\omega$, $\alpha$ range over the domain $\Im(\omega)>0$, $\alpha\in \C$.
Then as $v\to 0^+$, we have the estimate
\[
\log
\frac{\Gammat(e(v\alpha\omega);e(v\omega))}
     {\Gammar(\alpha)}
=
\frac{\pi\sqrt{-1}}{12 v\omega}
+
(\alpha-1/2)\log(2\pi v\omega/\sqrt{-1})
+
O(v),
\]
uniformly over compact subsets of the domain.
\end{cor}


\section{Generalized triangle inequalities}

When taking limits of elliptic hypergeometric integrals, the first step is
naturally to determine which part of the contour makes the most significant
contribution to the integral.  We first note the following consequence of
Proposition \ref{prop:ell_bound_hyper}.  Since in the sequel, we will only
be using gamma functions for $r\le 3$, we will write $\theta$ for
$\Gamma^{(1)}_e$, $\Gammae$ for $\Gamma^{(2)}_e$, and $\Gammae^+$ for
$\Gamma^{(3)}_e$, and similarly (with a subscript $h$) for the hyperbolic
versions; we will also omit the superscript $(r)$ on $P$, $Q$, and $R$.

\begin{cor}\label{cor:ell_bound_hyper}
For any parameters $\mu$, $\nu$, and any real number $x$, we have, as $v\to
0^+$, the estimates (uniform over compact subsets avoiding the poles)
\begin{align}
\frac{\Gammae(e(v\mu+x),e(v\nu-x);e(v\omega_1),e(v\omega_2))}
{e(R(v\mu;v\omega_1,v\omega_2)+R(v\nu;v\omega_1,v\omega_2))}
&=
O(e(\frac{\mu+\nu-\omega_1-\omega_2}{2v\omega_1\omega_2}\vartheta(x))),\notag\\
\frac{e(R(v\nu;v\omega_1,v\omega_2))
      \Gammae(e(v\mu+x);e(v\omega_1),e(v\omega_2))}
     {e(R(v\mu;v\omega_1,v\omega_2))
      \Gammae(e(v\nu+x);e(v\omega_1),e(v\omega_2))}
&=
O(e(\frac{\mu-\nu}{2v\omega_1\omega_2}\vartheta(x))),\notag\\
\frac{e(R(v\mu;v\omega_1,v\omega_2)+R(v\nu;v\omega_1,v\omega_2))}
{\Gammae(e(v\mu+x),e(v\nu-x);e(v\omega_1),e(v\omega_2))}
&=
O(e(\frac{\omega_1+\omega_2-\mu-\nu}{2v\omega_1\omega_2}\vartheta(x))),\notag
\end{align}
where $\vartheta(x)$ is the continuous, even, periodic function defined by
\[
\vartheta(x)=\{x\}(1-\{x\})=\{x\}\{-x\}.
\]
Similarly,
\[
e(-R(v\mu;v\omega))
\theta(e(v\mu+x);e(v\omega))
=
O(e(\vartheta(x)/2v\omega)),
\]
and likewise for the reciprocal (avoiding the poles).
\end{cor}

\begin{rem}
Note that in the above bounds, we can ignore $O(1)$ terms in $R$, and may
thus replace
\begin{align}
R(v\mu;v\omega) &\mapsto -1/12v\omega\notag\\
R(v\mu;v\omega_1,v\omega_2) &\mapsto
(\omega_1+\omega_2-2\mu)/24v\omega_1\omega_2.\notag
\end{align}
Also, although we assume $x$ real, the above estimates are clearly still
valid if we make a $O(v)$ perturbation to $x$ on the left-hand sides.
\end{rem}

We will thus require some inequalities involving this quantity $\vartheta$.

\begin{lem}\label{lem:gen_tri_ell}
For any sequences $c_1$,\dots,$c_n$, $d_1$,\dots,$d_n$ of real numbers,
we have the inequality
\[
 \sum_{1\le i,j\le n} \vartheta(c_i-d_j)
-\sum_{1\le i<j\le n} \vartheta(c_i-c_j)
-\sum_{1\le i<j\le n} \vartheta(d_i-d_j)
\ge
\vartheta(\sum_{1\le i\le n} c_i-d_i),
\]
with equality if and only if the sequences interlace in $\R/\Z$; that is,
iff they can be permuted so that either
\[
\{c_1\}\le \{d_1\}\le \{c_2\}\le\cdots\le\{d_{n-1}\}\le \{c_n\}\le \{d_n\}
\]
or
\[
\{d_1\}\le \{c_1\}\le \{d_2\}\le\cdots\le\{c_{n-1}\}\le \{d_n\}\le \{c_n\}.
\]
\end{lem}

\begin{proof}
First, observe that if two elements $c_i$, $d_j$ agree modulo $\Z$, then
their contributions to the inequality cancel, and the result thus follows
by induction.  We may therefore assume that $c_i\ne d_j$ for all $1\le
i,j\le n$.

Now, consider the asymptotics as $v\to 0^+$ of the case $\tau=\omega/2$ of
the determinant identity \cite{FrobeniusG:1882}:
\begin{align}
&\det_{1\le i,j\le n}(
\frac{e(-1/12v\omega)\theta(e(v\tau+c'_i-d'_j);e(v\omega))}
     {\theta(e(v\tau),e(c'_i-d'_j);e(v\omega))}
)\notag\\
&\qquad=
(-1)^{n(n-1)/2}
\frac{\theta(e(v\tau+\sum_i c'_i-d'_i);e(v\omega))}
     {e(n/12v\omega)\theta(e(v\tau);e(v\omega))}
\frac{\prod_{1\le i<j\le n} e(c'_j-d'_i)\theta(e(c'_i-c'_j),e(d'_i-d'_j);e(v\omega))}
{\prod_{1\le i,j\le n} \theta(e(c'_i-d'_j);e(v\omega))},
\notag
\end{align}
where real constants $(c'_i-c_i)/v$, $(d'_i-d_i)/v$ are chosen so that the
$2n$ quantities $c'_i$, $d'_i$ are all distinct for sufficiently small $v$.
Now, since $c_i\ne d_j$, we have
\[
\lim_{v\to 0^+}
\frac{e(-1/12v\omega)\theta(e(v\omega/2+c'_i-d'_j);e(v\omega))}
     {\theta(e(v\omega/2),e(c'_i-d'_j);e(v\omega))}
=
\frac{1}{2}\sgn(\{c_i\}-\{d_j\})e(-(\{c_i\}-\{d_j\})/2+1/4).
\]
In particular, the determinant converges to a well-defined limit.
Moreover, this limit is nonzero iff the sequences interlace, as follows by
considering the rescaled determinant $\det_{1\le i,j\le
  n}(\sgn(\{c_i\}-\{d_j\}))$.  (Indeed, if the sequences fail to interlace,
two rows or columns must agree; otherwise, the $n$ distinct rows are easily
verified to be linearly independent.)

On the other hand, we have the estimates
\begin{align}
\frac{\theta(e(v\omega/2+x);e(v\omega))}
     {\theta(e(v\omega/2);e(v\omega))}
&=
\Theta(e(\vartheta(x)/2v\omega))\notag\\
e(1/12v\omega)\theta(e(x);e(v\omega))
&=\Theta(e(\vartheta(x)/2v\omega))
\qquad \text{assuming}\ |\{x\}|,|\{1-x\}|=\Omega(v),
\notag
\end{align}
giving the estimate
\[
\Theta(e((\vartheta(\sum_i c_i-d_i)+\sum_{1\le i<j\le n}
(\vartheta(c_i-c_j)+\vartheta(d_i-d_j))-\sum_{1\le i,j\le n}\vartheta(c_i-d_j))/2v\omega))
\]
for the right-hand side.  Since $\Im(-1/\omega)>0$, this is bounded as
$v\to 0^+$ iff
\[
 \vartheta(\sum_i c_i-d_i)
+\sum_{1\le i<j\le n}(\vartheta(c_i-c_j)+\vartheta(d_i-d_j))
-\sum_{1\le i,j\le n}\vartheta(c_i-d_j)
\le 0,
\]
and is bounded away from 0 iff equality holds.  The result follows.
\end{proof}

\begin{rem}
More precise asymptotic calculations give the following well-known (and
easy) identity (valid for $n\ge 1$) as a limit of the elliptic Cauchy
determinant:
\[
\det_{1\le i,j\le n}(\sgn(x_i-y_j))
=
2^{n-1}
(-1)^{n(n-1)}
\prod_{1\le i,j\le n} \sgn(x_i-y_j)
\prod_{1\le i<j\le n} \sgn(x_i-x_j)\sgn(y_i-y_j)
\]
for interlacing sequences with distinct elements, and 0 otherwise.
Similarly, for more general values of $\tau$, we obtain the identity
\[
\det_{1\le i,j\le n}(x^{\sgn(c_i-d_j)})
=
(-1)^{n(n-1)/2} (x-1/x)^{n-1}
\sgn(\sum_i c_i-d_i) x^{\sgn(\sum_i c_i-d_i)}
\frac{\prod_{1\le i<j\le n} \sgn(c_i-c_j)\sgn(d_i-d_j)}
     {\prod_{1\le i,j\le n} \sgn(c_i-d_j)}.
\]
\end{rem}

We can also obtain a version with hyperoctahedral symmetry.

\begin{lem}\label{lem:gen_tri_ellC}
For any sequences $c_0$,\dots,$c_n$, $d_1$,\dots,$d_n$ of real numbers, we
have the inequality
\[
\sum_{0\le i\le n,1\le j\le n} \vartheta( c_i\pm d_j)
-\sum_{0\le i<j\le n} \vartheta(c_i\pm c_j)
-\sum_{1\le i<j\le n} \vartheta(d_i\pm d_j)
-\sum_{1\le i\le n} \vartheta(2d_i)
\ge
0,
\]
with equality iff the sequences can be permuted so that
\[
\min(\{\pm c_0\})
\le
\min(\{\pm d_1\})
\le
\min(\{\pm c_1\})
\le\cdots\le
\min(\{\pm c_{n-1}\})
\le
\min(\{\pm d_n\})
\le
\min(\{\pm c_n\}).
\]
Here $\vartheta(x\pm y):=\vartheta(x+y)+\vartheta(x-y)$ and $\min(\{\pm
x\}):=\min(\{x\},\{-x\})$.
\end{lem}

\begin{proof}
Apply the preceding lemma to the sequences $\pm c_i$ and $0,\pm d_i,1/2$,
and use the identity
\[
\vartheta(2x)=2(\vartheta(x)+\vartheta(x+1/2)-\vartheta(1/2)).
\]
The given equality condition simply restates the condition that $\pm c_i$
and $0,\pm d_i,1/2$ interlace.
\end{proof}

\begin{cor}\label{cor:gen_tri_ellC}
For any integer $n\ge 1$ and any sequence $e_1,\cdots,e_n$ of real numbers, we
have the inequality
\[
2\sum_{1\le i<j\le n} \vartheta(e_i\pm e_j)
-(n-1)\sum_{1\le i\le n} \vartheta(2e_i)
\ge
0
\]
with equality iff the sequence $\min(\{\pm e_i\})$ is constant.
\end{cor}

\begin{proof}
The case $c_0=0$, $n=1$ of Lemma \ref{lem:gen_tri_ellC} implies that
\[
\vartheta(x\pm y) + 2(\vartheta(x)-\vartheta(y))-\vartheta(2x)\ge 0
\] 
with equality iff $\min(\{\pm x\})\le \min(\{\pm y\})$.  Adding all
specializations of the form $(x,y)\mapsto (e_i,e_j)$ with $i\ne j$
gives the desired result.
\end{proof}

If we rescale $c_i,d_i\mapsto v c_i,vd_i$ and take $v\to 0^+$, the fact
that $\vartheta(vx)=v|x|-v^2 x^2$ for sufficiently small $v$ gives us the
following limit.

\begin{lem}\label{lem:gen_tri2}
For any sequences $c_1,\dots,c_n$, $d_1,\dots,d_n$, of real numbers, we
have the following inequality:
\[
\sum_{\substack{1\le i\le n\\1\le j\le n}} |c_i-d_j|
-\sum_{1\le i<j\le n} |c_i-c_j|
-\sum_{1\le i<j\le n} |d_i-d_j|
\ge 
|\sum_{1\le i\le n} c_i-\sum_{1\le i\le n} d_i|,
\]
with equality iff the sequences can be permuted so that
\[
c_1\le d_1\le\cdots\le c_n\le d_n
\]
or
\[
d_1\le c_1\le\cdots\le d_n\le c_n.
\]
\end{lem}

In particular, we have the following fact.

\begin{lem}\label{lem:gen_tri1}
For any sequences $c_0,\dots,c_n$, $d_1,\dots,d_n$, of real numbers, we
have the following inequality:
\[
\sum_{\substack{0\le i\le n\\1\le j\le n}} |c_i-d_j|
-\sum_{0\le i<j\le n} |c_i-c_j|
-\sum_{1\le i<j\le n} |d_i-d_j|
\ge 0,
\]
with equality iff the sequences can be permuted so that
\[
c_0\le d_1\le c_1\le\cdots\le c_{n-1}\le d_n\le c_n.
\]
\end{lem}

\begin{proof}
Choose a number $d_0$ such that $d_0<\min(c_0,\dots,c_n,d_1,\dots,d_n,\sum_i
c_i-\sum_{i>0} d_i)$, and apply Lemma \ref{lem:gen_tri2}.
\end{proof}

\begin{rem}
The case $n=1$ is of course just the usual triangle inequality in $\R$,
thus justifying the title of this section.
\end{rem}

%
%
%

\section{Hyperbolic limits}

Using the above asymptotic estimates for the hyperbolic and elliptic gamma
functions, we can obtain corresponding estimates for the various elliptic
hypergeometric integrals of \cite{Rains:Transformations} in the hyperbolic
limit.  In particular, in each case, it will turn out that up to an
explicit exponential factor, the elliptic integral converges exponentially
quickly to the hyperbolic integral.

Let us first consider the case of the Type I (perhaps better named
``elliptic Dixon'', see Corollary \ref{cor:Dixon} below and
\cite{DixonAL:1905}) integral with $BC_n$ symmetry, defined for all
nonnegative integers $m$, $n$, and parameters $p$, $q$, $u_0$\dots
$u_{2m+2n+3}$ satisfying
\[
0<|p|,|q|<1,\quad \prod_{0\le r\le 2m+2n+3} u_r = (pq)^{m+1}
\]
by the integral
\[
I^{(m)}_{BC_n}(u_0,u_1,\dots ;p,q)
:=
\frac{(p;p)^n(q;q)^n}{2^n n!}
\int_{C^n}
\frac{\prod_{1\le i\le n}\prod_{0\le r\le 2m+2n+3} \Gammae(u_r z_i^{\pm 1};p,q)}
{\prod_{1\le i<j\le n} \Gammae(z_i^{\pm 1}z_j^{\pm 1};p,q)
\prod_{1\le i\le n} \Gammae(z_i^{\pm 2};p,q)}
\prod_{1\le i\le n} \frac{dz_i}{2\pi\sqrt{-1}z_i},
\]
where the contour is chosen to contain all points of the form $p^i q^j
u_r$, $0\le i,j$, and exclude their reciprocals.

In the hyperbolic limit $p,q,u_r\to 1$, this gives rise to the following
limit.

\begin{thm}
Let $\mu_0$, $\mu_1$,\dots, $\mu_{2m+2n+3}$, $\omega_1$, $\omega_2$ be
parameters such that
\[
\Im(\omega_1),\Im(\omega_2)>0,\quad \sum_r \mu_r = (m+1)(\omega_1+\omega_2).
\]
Then as $v\to 0^+$,
\[
e(-2n\sum_r R(v\mu_r;v\omega_1,v\omega_2)+(2 n^2+n) R(0;v\omega_1,v\omega_2))
I^{(m)}_{BC_n}(e(v\mu_0),e(v\mu_1),\dots;e(v\omega_1),e(v\omega_2))
\]
converges uniformly exponentially (over compact subsets) to
\[
\frac{1}{(\sqrt{-\omega_1\omega_2})^n 2^n n!}
\int_{C^n}
\frac{\prod_{1\le i\le n}\prod_{0\le r\le 2m+2n+3} \Gammah(\mu_r\pm x_i;\omega_1,\omega_2)}
{\prod_{1\le i<j\le n} \Gammah(\pm x_i\pm x_j;\omega_1,\omega_2)
\prod_{1\le i\le n} \Gammah(\pm 2x_i;\omega_1,\omega_2)}
\prod_{1\le i\le n} dx_i,
\]
where the contour agrees with $\R$ outside a compact set, and is chosen to
contain all points of the form $\mu_r+i\omega_1+j\omega_2$, $i,j\ge 0$ and
exclude their negatives.
\end{thm}

\begin{proof}
We first observe that
\begin{align}
e(R(0;v\omega_1,v\omega_2))
(e(v\omega_1);e(v\omega_1))
(e(v\omega_2);e(v\omega_2))
&=
\lim_{z\to 0}
\frac{e(R(vz;v\omega_1,v\omega_2))}
{(1-e(vz))\Gammae(e(vz);e(v\omega_1),e(v\omega_2))}\notag\\
&\sim
\lim_{z\to 0}
\frac{1}{(1-e(vz))\Gammah(z;\omega_1,\omega_2)}\notag\\
&=
\frac{1}{v\sqrt{-\omega_1\omega_2}},\notag
\end{align}
with uniform exponentially small relative error as $v\to 0^+$.

For the remaining factors, we first assume that $\Im(\mu_r)>0$ for all $r$,
and thus the elliptic contour may be taken to be the unit circle.  Now, in
the elliptic integral, introduce the change of variables $z_i=e(x_i)$, and
thus $dz_i/2\pi\sqrt{-1}z_i = dx_i$; this replaces the unit circle by the
cube $[-1/2,1/2]^n$.  We next claim that if we restrict to the smaller cube
$[-1/4,1/4]^n$, the resulting error is uniformly exponentially small.
Indeed, we can use Corollary \ref{cor:ell_bound_hyper} to bound the
integrand on the full cube.  The $\mu$ factors satisfy
\[
e(-2R(v\mu_r;v\omega_1,v\omega_2))
\Gammae(e(v\mu_r\pm x_i);e(v\omega_1),e(v\omega_2))
=
O(
e(\frac{-\omega_1-\omega_2+2\mu_r}{2v\omega_1\omega_2} \vartheta(x_i))
)
\]
and thus, using the balancing condition,
\[
\prod_{0\le r\le 2m+2n+3}
e(-2R(v\mu_r;v\omega_1,v\omega_2))
\Gammae(e(v\mu_r\pm x_i);e(v\omega_1),e(v\omega_2))
=
O(
e(
\frac{-\omega_1-\omega_2}{2v\omega_1\omega_2} 2(n+1)\vartheta(x_i)
)
).
\]
Similarly, the remaining univariate factors satisfy
\[
e(2R(0;v\omega_1,v\omega_2))
\Gammae(e(\pm 2x_i);e(v\omega_1),e(v\omega_2))^{-1}
=
O(
e(\frac{-\omega_1-\omega_2}{2v\omega_1\omega_2} (-\vartheta(2x_i)))
)
\]
and, for $i<j$, the cross factors satisfy
\[
e(4R(0;v\omega_1,v\omega_2))
\Gammae(e(\pm x_i\pm x_j);e(v\omega_1),e(v\omega_2))^{-1}
=
O(
e(\frac{-\omega_1-\omega_2}{2v\omega_1\omega_2} (-\vartheta(x_i\pm x_j)))
)
\]
Combining these bounds, we find that the integrand is uniformly
\[
O(
e(
\frac{-\omega_1-\omega_2}{2v\omega_1\omega_2}
\bigl(
\sum_{1\le i\le n} (2n+2)\vartheta(x_i)
-\sum_{1\le i\le n} \vartheta(2x_i)-\sum_{1\le i<j\le n} \vartheta(x_i\pm
x_j)
\bigr)
))
\]
Since
\[
\Im(\frac{-\omega_1-\omega_2}{\omega_1\omega_2})=\Im(-1/\omega_1)+\Im(-1/\omega_2)>0,
\]
the bound is maximized when
\[
\sum_{1\le i\le n} (2n+2)\vartheta(x_i)
-\sum_{1\le i\le n} \vartheta(2x_i)-\sum_{1\le i<j\le n} (\vartheta(x_i+
x_j)+\vartheta(x_i-x_j))
\]
is minimized, which in turn happens when $x_1=x_2=\cdots=x_n=0$, by
Lemma \ref{lem:gen_tri_ellC} applied to the case $c_i\equiv 0$,
$d_i=x_i$.  In particular, the integrand is exponentially small
everywhere else, and thus restricting to $|x_i|\le 1/4$ introduces an
exponentially small error.

At this point, using Proposition \ref{prop:ell_to_hyper2} allows us to
replace the gamma functions in the integrand with hyperbolic gamma
functions (times an exponential factor that turns out to be trivial).  The
factor $v^{-n}$ from $((p;p)(q;q))^n$ can be absorbed in rescaling the
variables of integration; we thus obtain the restriction of the desired
integral to the cube $[-1/4v,1/4v]$.  But again we can
bound the integrand, this time using Corollary \ref{cor:hyperb_asympt}, and
find the uniform bound
\[
O(
e(
\frac{-\omega_1-\omega_2}{2\omega_1\omega_2}
\left(
\sum_{1\le i\le n} (2n+2)|x_i|
-\sum_{1\le i\le n} |2x_i|-\sum_{1\le i<j\le n} (|x_i+x_j|+|x_i-x_j|)
\right)
)
),
\]
so the omitted tail is again uniformly exponentially small.

For the general case, we note that if $C$ is a valid choice of contour for
the hyperbolic integral, then for sufficiently small $v$, the image of the
subcontour $[-1/2v,1/2v]$ under $x\mapsto e(vx)$ is a valid choice of
contour for the elliptic integral.  The result agrees with the unit circle
outside a neighborhood of size $O(v)$ of 1; as a result, the difference
from the unit circle has no effect on the asymptotics.
\end{proof}

If we denote the above integral by $I^{(m)}_{BC_n;h}$, we have the
following corollary, obtained as the limit of the corresponding identity
for the elliptic case; note that we do not need to compare the exponential
factors on both sides, since both sides must agree throughout and have
generically nonzero limits.

\begin{cor}
Let $\mu_0$, $\mu_1$,\dots, $\mu_{2m+2n+3}$, $\omega_1$, $\omega_2$ be
parameters such that
\[
\Im(\omega_1),\Im(\omega_2)>0,\quad \sum_r \mu_r = (m+1)(\omega_1+\omega_2).
\]
Then
\[
I^{(m)}_{BC_n;h}(\dots,\mu_r,\dots;\omega_1,\omega_2)
=
\prod_{0\le r<s\le 2m+2n+3} \Gammah(\mu_r+\mu_s;\omega_1,\omega_2)
I^{(n)}_{BC_m;h}(\dots,\frac{\omega_1+\omega_2}{2}-\mu_r,\dots;\omega_1,\omega_2),
\]
and in particular
\[
I^{(0)}_{BC_n;h}(\dots,\mu_r,\dots;\omega_1,\omega_2)
=
\prod_{0\le r<s\le 2n+3} \Gammah(\mu_r+\mu_s;\omega_1,\omega_2).
\]
\end{cor}

\begin{rem}
As van Diejen and Spiridonov \cite{vanDiejenJF/SpiridonovVP:2005} observed
for the Type II evaluation, one can also prove hyperbolic results by simply
replacing the arguments of \cite{Rains:Transformations} by appropriate
limits, rather than taking limits directly.  Those arguments depend
strongly on the fact that the set $p^\Z q^\Z$ generically has finite limit
points (in fact, is generically dense), which makes analytic continuation
trivial.  In the hyperbolic setting, the corresponding set
$\Z\omega_1+\Z\omega_2$ is never dense, and only has a limit point when
$\omega_1/\omega_2$ is real irrational, so an additional analytic
continuation argument is needed to extend to generic moduli.
\end{rem}

A similar argument will work in the other cases; some technical issues do
arise, however, so it is worth discussing those cases as well.

For the Type II (again, the name ``elliptic Selberg'' might be better)
integral, the main complication is that without an additional condition on
the parameters, the integrand is not maximized near $z_i\equiv 1$.  We have
the following result.  Define a family of integrals
\begin{align}
\II^{(m)}_{BC_n}&(u_0,u_1,\dots,u_{2m+5};t;p,q)\notag\\
&:=
\frac{(p;p)^n(q;q)^n\Gammae(t;p,q)^n}{2^n n!}
\int_{C^n}
\prod_{1\le i<j\le n} \frac{\Gammae(t z_i^{\pm 1} z_j^{\pm
    1};p,q)}{\Gammae(z_i^{\pm 1} z_j^{\pm 1};p,q)}
\prod_{1\le i\le n}
\frac{\prod_{0\le r\le 2m+5} \Gammae(u_r z_i^{\pm 1};p,q)}{\Gammae(z_i^{\pm
    2};p,q)}
\frac{dz_i}{2\pi\sqrt{-1}z_i},
\notag
\end{align}
on the domain $t^{2n-2}\prod_r u_r = (pq)^{m+1}$, $0<|p|,|q|,|t|<1$, where
the contour $C$ satisfies $C=C^{-1}$, and for all $i,j\ge 0$, contains the
points $p^i q^j u_r$ as well as the contour $p^i q^j t C$.

\begin{thm}
Let $\mu_0$, $\mu_1$,\dots, $\mu_{2m+5}$, $\tau$, $\omega_1$,
$\omega_2$ be parameters such that
\[
\Im(\tau),\Im(\omega_1),\Im(\omega_2)>0,\quad (2n-2)\tau+\sum_r \mu_r = (m+1)(\omega_1+\omega_2),
\]
and satisfying the convergence condition
\[
\Im(\frac{-(n-1)\tau-\omega_1-\omega_2}{\omega_1\omega_2})>0.
\]
Then as $v\to 0^+$,
\begin{align}
e(-2n\sum_r R(v\mu_r;v\omega_1,v\omega_2)+2n^2R(0;v\omega_1,v\omega_2)
  -2n(n-1)R(v\tau;v\omega_1,v\omega_2))\qquad&\notag\\
\cdot\II^{(m)}_{BC_n}(e(v\mu_0),e(v\mu_1),\dots;e(v\tau);e(v\omega_1),e(v\omega_2))&
\notag
\end{align}
converges uniformly exponentially (over compact subsets) to
\[
\frac{\Gammah(\tau;\omega_1,\omega_2)^n}{(\sqrt{-\omega_1\omega_2})^n 2^n n!}
\int_{C^n}
\prod_{1\le i<j\le n}
\frac{\Gammah(\tau\pm x_i\pm x_j;\omega_1,\omega_2)}
     {\Gammah(\pm x_i\pm x_j;\omega_1,\omega_2)}
\prod_{1\le i\le n}
\frac{\prod_{0\le r\le 2m+5} \Gammah(\mu_r\pm x_i;\omega_1,\omega_2)}
     {\Gammah(\pm 2x_i;\omega_1,\omega_2)}
dx_i,
\]
where
the contour $C=-C$ agrees with $\R$ outside a compact set and for all
$i,j\ge 0$ contains the points $i\omega_1+j\omega_2+\mu_r$ as well as the
contour $i\omega_1+j\omega_2+\tau+C$.
\end{thm}

\begin{proof}
Again we change variables to $z_i=e(x_i)$ and integrate over the cube
$[-1/2,1/2]^n$; we may also freely assume $n>1$, as the case $n=1$ has
already been dealt with.  In this case, we find that the integrand is
uniformly bounded by
\begin{align}
O(e(&
\frac{-\omega_1-\omega_2}{2(n-1)v\omega_1\omega_2}
\Bigl(
2\sum_{1\le i<j\le n} \vartheta(x_i\pm x_j)
-(n-1)\sum_{1\le i\le n}\vartheta(2x_i)
\Bigr)\notag\\
{}+{}&
\frac{-(n-1)\tau-\omega_1-\omega_2}{(n-1)v\omega_1\omega_2}
\Bigl(
2(n-1)\sum_{1\le i\le n}\vartheta(x_i)
-
\sum_{1\le i<j\le n} \vartheta(x_i\pm x_j)
\Bigr)
)).
\notag
\end{align}
By Corollary \ref{cor:gen_tri_ellC}, the first $\vartheta$ sum is $\ge 0$,
with equality iff the sequence $|x_i|$ is constant; by the case $d_i\equiv
0$ of Lemma \ref{lem:gen_tri_ellC}, the second $\vartheta$ sum is $\ge 0$,
with equality iff at most one of the $x_i$ is nonzero.  It follows that the
integrand is exponentially small unless both conditions are satisfied;
i.e., unless $x_i\equiv 0$.  The remainder of the proof is as above.
\end{proof}

\begin{rem}
Note that it also follows from the above proof that the convergence
condition is necessary for the integrand to be localized.  One can readily
arrange for equality to hold in the first sum, but not the second, at which
point if $\Im(-((n-1)\tau+\omega_1+\omega_2)/\omega_1\omega)<0$, the
integrand is exponentially larger than its value near $x_i\equiv 0$.
\end{rem}

Let $\II^{(m)}_{BC_n;h}(\mu_0,\dots;\tau;\omega_1,\omega_2)$ denote the
above hyperbolic integral, as a meromorphic function on the domain
\[
\Im(\tau),\Im(\omega_1),\Im(\omega_2),\Im(\frac{-(n-1)\tau-\omega_1-\omega_2}{\omega_1\omega_2})>0,\qquad
(2n-2)\tau+\sum_r \mu_r = (m+1)(\omega_1+\omega_2).
\]

\begin{cor}
Let $\mu_0$, $\mu_1$,\dots, $\mu_5$, $\tau$, $\omega_1$, $\omega_2$ be
parameters such that
\[
\Im(\tau),\Im(\omega_1),\Im(\omega_2)>0,\quad
(2n-2)\tau+\sum_r \mu_r = \omega_1+\omega_2,
\]
and satisfying the convergence condition
\[
\Im(\frac{-(n-1)\tau-\omega_1-\omega_2}{\omega_1\omega_2})>0.
\]
Then
\[
\II^{(0)}_{BC_n;h}(\mu_0,\dots,\mu_5;\tau;\omega_1,\omega_2)
=
\prod_{0\le i<n} \Gammah((i+1)\tau;\omega_1,\omega_2)
\prod_{0\le r<s\le 5} \Gammah(i\tau+\mu_r+\mu_s;\omega_1,\omega_2).
\]
\end{cor}

\begin{cor}
For parameters $\mu_0,\dots,\mu_7$, $\tau$, $\omega_1$, $\omega_2$ such
that there exists an integer $n$ (necessarily unique) with
\[
(2n+2)\tau+\sum_r \mu_r = 2(\omega_1+\omega_2)
\]
and
\[
\Im(\tau),\Im(\omega_1),\Im(\omega_2),\Im(\frac{-(n-1)\tau-\omega_1-\omega_2}{\omega_1\omega_2})>0,
\]
define
\[
\tilde{\II}_h(\mu_0,\dots,\mu_7;\tau;\omega_1,\omega_2)
:=
\bigl(\prod_{0\le r<s\le 7}
\Gammah^+(\tau+\mu_r+\mu_s;\tau,\omega_1,\omega_2)
\bigr)
\II^{(1)}_h(\tau/2+\mu_0,\dots,\tau/2+\mu_7;\tau;\omega_1,\omega_2).
\]
Then $\tilde{\II}_h$ is invariant under the natural action of the Weyl
group $E_7$; in other words, it satisfies the identities
\[
\tilde{\II}_h(\mu_0,\dots,\mu_7;\tau;\omega_1,\omega_2)
=
\tilde{\II}_h(\mu_0+\nu,\dots,\mu_3+\nu,\mu_4-\nu,\dots,\mu_7-\nu;\tau;\omega_1,\omega_2)
\]
where $\nu=(\mu_4+\mu_5+\mu_6+\mu_7-\mu_0-\mu_1-\mu_2-\mu_3)/4$;
\[
\tilde{\II}_h(\mu_0,\dots,\mu_7;\tau;\omega_1,\omega_2)
=
\tilde{\II}_h(\nu-\mu_0,\dots,\nu-\mu_3,\nu'-\mu_4,\dots,\nu'-\mu_7;\tau;\omega_1,\omega_2),
\]
where $\nu=(\mu_0+\mu_1+\mu_2+\mu_3)/2$,
$\nu'=(\mu_4+\mu_5+\mu_6+\mu_7)/2$; and
\[
\tilde{\II}_h(\mu_0,\dots,\mu_7;\tau;\omega_1,\omega_2)
=
\tilde{\II}_h(\nu-\mu_0,\dots,\nu-\mu_7;\tau;\omega_1,\omega_2),
\]
where $\nu=(\mu_0+\mu_1+\dots+\mu_7)/2$; as well as invariance under
permutations of $\mu_0$ through $\mu_7$.
\end{cor}

\begin{rem}
Similarly, the other double coset of $E_7$ in $E_8$ that gives rise to
(dimension-altering) transformations of the elliptic integral also gives
rise to transformations of the hyperbolic integral; we omit the obvious
details.  The key observation is that the overall exponential factor that
arises when taking the limit is, once one solves for $n$, a function of
$\sum_i \mu_i^2$, and is thus $E_8$-invariant.  The work of
\cite{Rains:Recurrences} on recurrences also descends to the hyperbolic
case; in particular, for $\tau=\omega_2$, one obtains a tau-function for a
hyperbolic analogue of the elliptic Painlev\'e equation.
\end{rem}

For the $A_n$ integral, the difficulty is that the elliptic integral has a
condition $\prod_i z_i = 1$, which in $x_i$ coordinates, becomes $\sum_i
x_i\in \Z$; this introduces extra complications when maximizing the
integrand.  Recall that the $A_n$ integral is defined by
\begin{align}
I^{(m)}_{A_n}(u_0,\dots u_{m+n+1};&v_0,\dots v_{m+n+1};p,q)\notag\\
&:=
\frac{(p;p)^n(q;q)^n}{(n+1)!}
\int_{\prod_{0\le i\le n} z_i=1}
\frac{
\prod_{0\le i\le n}
\prod_{0\le r<m+n+2} \Gammae(u_r z_i,v_r/z_i;p,q)
}{
\prod_{0\le i<j\le n} \Gammae(z_i/z_j,z_j/z_i;p,q)
}
\prod_{1\le i\le n} \frac{dz_i}{2\pi\sqrt{-1}z_i},\notag
\end{align}
for $0<|p|,|q|<1$, $0<|u_0|,\dots,|u_{m+n+1}|,|v_0|,\dots,|v_{m+n+1}|<1$,
$\prod_i u_iv_i=(pq)^{m+1}$.  (It follows from general principles that this
can be extended to a meromorphic function on the domain $0<|p|,|q|<1$,
$\prod_i u_iv_i=(pq)^{m+1}$, but the condition on the contour is
complicated.)

\begin{thm}
Let $\mu_0$, $\mu_1$,\dots, $\mu_{m+n+1}$, $\nu_0$, $\nu_1$,\dots,
$\nu_{m+n+1}$, $\omega_1$, $\omega_2$ be parameters in the upper
half-plane such that
\[
\sum_r \mu_r+\nu_r = (m+1)(\omega_1+\omega_2)
\]
Then as $v\to 0^+$,
\begin{align}
e(-(n+1)\sum_r (R(v\mu_r;v\omega_1,v\omega_2)+R(v\nu_r;v\omega_1,v\omega_2))
  +(n^2+2n)R(0;v\omega_1,v\omega_2))\qquad&\notag\\
\cdot I^{(m)}_{A_n}(e(v\mu_0),e(v\mu_1),\dots;e(v\nu_0),e(v\nu_1),\dots;e(v\omega_1),e(v\omega_2))&
\notag
\end{align}
converges uniformly exponentially (over compact subsets) to
\[
\frac{1}{(\sqrt{-\omega_1\omega_2})^n(n+1)!}
\int_{\sum_{0\le i\le n} x_i = 0}
\frac{
\prod_{0\le i\le n}
\prod_{0\le r<m+n+2} \Gammah(\mu_r+x_i,\nu_r-x_i;\omega_1,\omega_2)
}{
\prod_{0\le i<j\le n} \Gammah(x_i-x_j,x_j-x_i;\omega_1,\omega_2)
}
\prod_{1\le i\le n} dx_i.
\]
\end{thm}

\begin{proof}
If we perform the change of variables $z_i=e(x_i)$ in the elliptic
integral, the result is an integral over the domain
\[
-1/2\le x_0,x_1,x_2,\dots,x_n\le 1/2;\quad \sum_{0\le i\le n} x_i\in Z,
\]
a disjoint union of polytopes.  Over the entire cube, we find that the
integrand is uniformly
\[
O(e(
\frac{-\omega_1-\omega_2}{2v\omega_1\omega_2}
\bigl(
\sum_{0\le i\le n} (n+1)\vartheta(x_i)
-
\sum_{0\le i<j\le n} \vartheta(x_i-x_j)
\bigr)
)).
\]
Now, we find from the case $d_i\equiv 0$ of Lemma \ref{lem:gen_tri_ell}
that
\[
\sum_{0\le i\le n} (n+1)\vartheta(x_i)
-
\sum_{0\le i<j\le n} \vartheta(x_i-x_j)
\ge
\vartheta(\sum_i x_i)
\ge
0
\]
with equality iff $x_0,\dots,x_n$ interlaces with $0,\dots,0$ and $\sum_i
x_i\in \Z$; i.e., iff $x_i\equiv 0$.  We thus conclude that the integral
over the polytope
\[
-1/4\le x_0,x_1,\dots,x_n\le 1/4;\quad \sum_i x_i = 0
\]
is uniformly exponentially close to the original integral.
Thus, as above, the theorem reduces to showing that the hyperbolic integral
decays exponentially.  This in turn reduces to the identity
\[
(n+1)\sum_{0\le i\le n} |x_i| - \sum_{0\le i<j\le n} |x_i-x_j|
\ge 0
\]
with equality only when $x_1=x_2=\cdots x_n=0$.
\end{proof}

The remaining issue in degenerating \cite{Rains:Transformations} to the
hyperbolic level is the degeneration of the biorthogonal functions
constructed there.  The primary difficulty is that the construction of
those functions in \cite{Rains:Transformations} does not give rise to good
uniform asymptotics.  However, we can still establish the following.

\begin{thm}\label{thm:biorth_lim}
Let the parameters $\tau_0$, $\tau_1$, $\tau_2$, $\tau_3$, $\mu_0$,
$\mu_1$, $\tau$, $\omega_1$, $\omega_2$ be parameters with $\tau$,
$\omega_1$, $\omega_2$ in the upper half-plane such that
\[
(2n-2)\tau+\tau_0+\tau_1+\tau_2+\tau_3+\mu_0+\mu_1 = \omega_1+\omega_2.
\]
Then for any partition pair $\blambda$, and for generic values of the
parameters, the biorthogonal function
\[
\tilde{\cal
  R}^{(n)}_{\blambda}(\dots,e(x_i),\dots;e(v\tau_0){:}e(v\tau_1),e(v\tau_2),e(v\tau_3);e(v\mu_0),e(v\mu_1);e(v\tau);e(v\omega_1),e(v\omega_2))
\]
is uniformly bounded for $(x_1,\dots,x_n)\in D(v)^n$, where $D(v)$ is a
region of the form $-1-vC_1\le \Re(x)\le vC_1$, $|\Im(x)|\le vC_2$, and
excluding a hole of radius $\epsilon v$ around every pole of the
biorthogonal function.  Moreover, there exists a function
\[
\tilde{\cal
  R}^{(n)}_{\blambda;h}(\dots,x_i,\dots;\tau_0{:}\tau_1,\tau_2,\tau_3;\mu_0,\mu_1;\tau;\omega_1,\omega_2)
\]
such that as $v\to 0^+$,
\begin{align}
&\tilde{\cal
  R}^{(n)}_{\blambda}(\dots,e(x_i),\dots;e(v\tau_0){:}e(v\tau_1),e(v\tau_2),e(v\tau_3);e(v\mu_0),e(v\mu_1);e(v\tau);e(v\omega_1),e(v\omega_2))\notag\\
{}-{}&
\tilde{\cal
  R}^{(n)}_{\blambda;h}(\dots,x_i/v,\dots;\tau_0{:}\tau_1,\tau_2,\tau_3;\mu_0,\mu_1;\tau;\omega_1,\omega_2)
\notag
\end{align}
converges exponentially to 0, uniformly for $x$ in a compact subset of the
domain $-1<\Re(x)<1$.
\end{thm}

\begin{proof}
We first observe that the claims of the theorem are certainly true if we
replace $\tilde{\cal R}^{(n)}_{\blambda}$ by a product of functions of the
form
\[
e(2n(R(v\beta;v\omega_1)-R(v\alpha;v\omega_1)))\prod_{1\le i\le n}
\frac{\theta(e(v\alpha) z_i^{\pm 1};e(v\omega_1))}{\theta(e(v\beta) z_i^{\pm 1};e(v\omega_1))},
\]
or similarly for $\omega_2$.  In particular, it was established in
\cite{Rains:Transformations} that there exist functions
$F^{(n)}_{\blambda}(\mu_0{:};\tau;\omega_1,\omega_2)$ of the
above form such that there exists an expansion
\[
\tilde{\cal R}^{(n)}_{\blambda} = \sum_{\bmu\subset\blambda}
C_{\blambda\bmu} F_{\blambda}
\]
for some coefficients $C_{\blambda\bmu}$ independent of $z_i$.  It thus
remains only to show that for generic parameters, these coefficients
$C_{\blambda\bmu}$ converge exponentially.  Moreover, the action of the
integral operators of \cite{Rains:Transformations} can be computed
explicitly in the $F_{\blambda}$ basis, and the coefficients of the
corresponding matrices converge exponentially (to a triangular matrix with
generically nonzero diagonal).  Thus the generalized eigenvalue equations
satisfied by $\tilde{\cal R}^{(n)}_{\blambda}$ set up linear equations in
the $C_{\blambda\bmu}$ with exponentially converging coefficients.  Since
the limits of the generalized eigenvalues are generically distinct, the
limiting linear equations are generically nonsingular, and the result
follows.
\end{proof}

\begin{rem}
The unviariate hyperbolic biorthogonal function
\[
\tilde{\cal R}^{(1)}_{\blambda;h}(x;\tau_0{:}\tau_1,\tau_2,\tau_3;\mu_0,\mu_1;\tau;\omega_1,\omega_2)
\]
was discussed in \cite[\S 8.3]{SpiridonovVP:2005}.
\end{rem}

Note in particular that if we multiply the integrand of either $BC_n$
integral by a function satisfying such convergence properties, the
resulting integral will also converge exponentially (assuming the unadorned
integral so converges).  Also, a similar argument works for the
interpolation functions (which as special cases of the biorthogonal
functions do not quite fall under the above generic result, but again
satisfy suitable integral equations).  As a result, every identity of
\cite{Rains:Transformations} involving such functions converges
exponentially (possibly with an explicit factor of the form $\exp(a+bv)$)
to a corresponding hyperbolic limit.  One should note (as observed in
\cite[\S 8.3]{SpiridonovVP:2005}) that further degeneration of the
parameters can lead to convergence issues, as without the moderating effect
of the poles, the biorthogonal functions grow exponentially as
$|\Re(x)|\to\infty$.

\section{Trigonometric limits}

The main difficulty with the trigonometric limit $p\to 0$ is that the
general case of the transformations involves parameters tending to
infinity, making the contour ill-behaved in the limit.  This can be fixed
at the expense of breaking the symmetry of the integrand.

Recall that for the type $I$ $BC_n$ integral, the parameters are
constrained to satisfy the balancing condition
\[
\prod_{0\le r\le 2n+2m+3} u_r = (pq)^{m+1}.
\]
The natural way to satisfy this in the $p\to 0$ limit is for $2n+m+3$ of
the parameters to be $\Theta(1)$, while the remaining $m+1$ parameters are
$\Theta(p)$.  This then makes the $p\to 0$ limit of the integral trivial to
compute.

\begin{thm}\label{thm:triglim_IC1}
For any parameters $u_0,\dots,u_{2n+m+2}$, $v_0,\dots,v_m$, $q$ satisfying
\[
|q|<1,\quad \prod_{0\le r\le 2n+m+2} u_r = \prod_{0\le r\le m} v_r,
\]
we have the limit
\begin{align}
\lim_{p\to 0}
I^{(m)}_{BC_n}&(u_0,\dots,u_{2n+m+2},pq/v_0,\dots,pq/v_m;p,q)\notag\\
&=
\frac{(q;q)^n}{2^n n!}
\int_{C^n}
\prod_{1\le i<j\le n} \Gammat(z_i^{\pm 1}z_j^{\pm 1};q)^{-1}
\prod_{1\le i\le n}
\frac{ \prod_{0\le r\le 2n+m+2} \Gammat(u_r z_i^{\pm 1};q)}
     {\Gammat(z_i^{\pm 2};q) \prod_{0\le r\le m} \Gammat(v_r z_i^{\pm 1};q)}
\frac{dz_i}{2\pi\sqrt{-1}z_i},\notag
\end{align}
where the contour contains all points of the form $p^iq^ju_r$, $i,j\ge 0$,
and excludes their reciprocals.
\end{thm}

\begin{proof}
This follows immediately from the facts that as $p\to 0$,
\begin{align}
\Gammae(x;p,q)^{\pm 1} &= \Gammat(x;q)^{\pm 1}(1+O(p))\notag\\
\Gammae(pq/x;p,q)^{\pm 1} &= \Gammat(x;q)^{\mp 1}(1+O(p)),\notag
\end{align}
uniformly in $x$ away from the poles.
\end{proof}

Unfortunately, the right-hand side of the type I transformation involves
parameters
\[
(pq)^{1/2}/u_0,\dots,(pq)^{1/2}/u_{2n+m+2},(pq)^{-1/2}v_0,\dots,(pq)^{-1/2}v_m,
\]
which as mentioned above gives an apparently ill-behaved limit.  The
primary difficulty is that the divergent parameters not only deform the
contour, but in fact pinch the contour in the limit, making it approach
both 0 and infinity.  It turns out, however, that there is a way to break
the symmetry in such a way as to eliminate half of the offending poles,
thus allowing the contour to be renormalized, giving a well-behaved limit.

The key fact is the following identity of $q$-elliptic functions.  Here
$R(z_i)$ denotes the operator such that $R(z_i)f(z_i)=f(1/z_i)$.

\begin{lem}\label{lem:diffop_ICn}
For any parameters $u_0,\dots,u_{n+1}$, $q$ we have the identity
\[
\prod_{1\le i\le n} (1+R(z_i))
\frac{
\theta(\prod_{0\le r\le n+1} u_r/\prod_{1\le i\le n} z_i;q)
\prod_{1\le i\le n} \prod_{0\le r\le n+1} \theta(u_r z_i;q)
}
{\prod_{1\le i\le j\le n} \theta(z_iz_j;q)}
=
\prod_{0\le r<s\le n+1} \theta(u_ru_s;q).
\]
\end{lem}

\begin{proof}
The left-hand side can be expressed as a sum of $2^n$ terms, all of which
are elliptic functions in $z$ with respect to multiplication by $q$, and
thus the sum is also an elliptic function.  Moreover, since the original
function is invariant under permutations, the sum is invariant under the
action of $BC_n$.  In particular, the order of the sum along each
reflection hyperplane must be even; since the summands have at most simple
poles there, it follows that the sum is constant.  The constant can be
recovered by taking $z_i = u_i$, making all but one summand vanish.
\end{proof}

\begin{lem}\label{lem:nonsym_ICn}
For any nonzero parameters $t_0,\dots,t_n$, $u_0,\dots,u_{n+m+1}$,
$v_0,\dots,v_m$, $p$, $q$ with
\[
0<|p|,|q|<1,\quad
\prod_{0\le r\le n} t_r = \prod_{0\le r\le n+m+1}u_r\prod_{0\le r\le m}v_r
\]
and any complex parameters $a,w\ne 0$, we have the
identity
\begin{align}
\prod_{0\le r<s\le n} \theta(t_rt_s/a;q)
I^{(m)}_{BC_n}&(t_0/a^{1/2},\dots,t_n/a^{1/2},
                a^{1/2}/u_0,\dots,a^{1/2}/u_{m+n+1},
                pq/a^{1/2}v_0,\dots,pq/a^{1/2}v_m
;p,q)\notag\\
=
\frac{(p;p)^n(q;q)^n}{n!}
\int_{C^n}&
\frac{\theta(\prod_{0\le r\le n} t_r/w\prod_{1\le i\le n} z_i;q)
      \prod_{1\le i\le n} \theta(z_i/w;q)}
     {\prod_{0\le r\le n} \theta(t_r/w;q)}
\frac{\prod_{1\le i\le j\le n} \theta(pz_iz_j/a;p)}
{\prod_{1\le i<j\le n} \Gammae((z_i/z_j)^{\pm 1};p,q)}\notag\\
&\prod_{1\le i\le n}
  \prod_{0\le r\le n} \Gammae(p t_rz_i/a,t_r/z_i;p,q)
  \prod_{0\le r\le m+n+1} \frac{\Gammae(z_i/u_r;p,q)}{\Gammae(pqz_iu_r/a;p,q)}
\notag\\
&\phantom{\prod_{1\le i\le n}}
  \prod_{0\le r\le m} \frac{\Gammae(pq z_i/av_r;p,q)}{\Gammae(z_iv_r;p,q)}
  \frac{dz_i}{2\pi\sqrt{-1}z_i}
\notag
\end{align}
where the contour contains all points of the form $p^iq^jt_r$,
$p^iq^ja/u_r$, $p^{i+1}q^{j+1}/v_r$, $i,j\ge 0$, and excludes all points of
the form $p^{-1-i}q^{-j}a/t_r$, $p^{-i}q^{-j}u_r$, $p^{-i-1}q^{-j-1}av_r$,
$i,j\ge 0$.
\end{lem}

\begin{proof}
If we multiply the integrand on the left-hand side by the case
$u_r = a^{-1/2} v_r$, $0\le r\le n$, $u_{n+1}=a^{1/2}/w$ of the
lemma, the symmetry of the integrand implies that $\prod_{1\le i\le
  n}(1+R(z_i))$ can be replaced by $2^n$.  Shifting the variables of
integration by $z_i\to a^{-1/2}z_i$ gives the right-hand side, up to a
shift in contour with no effect on the integral.
\end{proof}

\begin{rems}
Note that for specific choices of $w$, the contour condition may
conceivably be weakened; the point is that the $w$-dependent factors can
cancel out poles of the integrand, making the corresponding constraints on
the contour superfluous.  In particular, for certain specializations of the
parameters, it can be the case that the contour conditions for generic $w$
are inconsistent, but a suitable choice of $w$ makes the integral
well-defined.
\end{rems}

\begin{rems}
If we multiply the integrand on the left by a symmetric function $f$
(adjusting the contour conditions accordingly), the effect is to multiply
the nonsymmetric integrand by $f(\dots a^{-1/2}z_i\dots)$, with suitable
contour conditions.
\end{rems}

This makes the limit $p\to 0$ trivial again, as long as $|pq|\le |a|<1$.
Taking $a=pq$ gives the following.

\begin{thm}\label{thm:triglim_IC2}
For any nonzero parameters $u_0,\dots,u_{2n+m+2}$, $v_0,\dots,v_m$, $q$ with
\[
0<|q|<1,\quad \prod_{0\le r\le 2n+m+2}u_r = \prod_{0\le r\le m} v_r,
\]
we have the limit
\begin{align}
\lim_{p\to 0}
\prod_{0\le r<s\le m} \theta(v_rv_s/pq;q)
I^{(n)}_{BC_m}&((pq)^{1/2}/u_0,\dots,(pq)^{1/2}/u_{2n+m+2},(pq)^{-1/2}v_0,\dots,(pq)^{-1/2}v_m;p,q)\notag\\
{}=
\frac{(q;q)^m}{m!}
\int_{C^m}&
\frac{
\theta(\prod_{0\le r\le m} v_r/w\prod_{1\le i\le m} z_i;q)
\prod_{1\le i\le m} \theta(z_i/w;q)
}
{\prod_{0\le r\le m} \theta(v_r/w;q)}
\frac{\prod_{1\le i\le j\le m} (1-z_iz_j/q)}
{\prod_{1\le i<j\le m} \Gammat((z_i/z_j)^{\pm 1};q)}\notag\\
&\prod_{1\le i\le m}
  \prod_{0\le r\le 2n+m+2} \frac{\Gammat(z_i/u_r;q)}{\Gammat(z_iu_r;q)}
  \prod_{0\le r\le m} \Gammat(v_r z_i/q,v_r/z_i;q)
  \frac{dz_i}{2\pi\sqrt{-1}z_i},
\notag
\end{align}
where the contour contains all points of the form
$q^j v_r$, $j\ge 0$, and excludes all points of the form
$q^{-j}u_r$, $q^{1-j}/v_r$, $j\ge 0$.
\end{thm}

\begin{cor}\label{cor:triglim_IC12}
The trigonometric integral of Theorem \ref{thm:triglim_IC1} is equal to 
\[
\prod_{0\le r<s\le 2n+m+2} \Gammat(u_ru_s;q)
\prod_{\substack{0\le r\le 2n+m+2\\0\le s\le m}} \Gammat(v_s/u_r;q)^{-1}
\prod_{0\le r<s\le m} \Gammat(v_rv_s/q;q)^{-1}
\]
times the trigonometric integral of Theorem \ref{thm:triglim_IC2}.
\end{cor}

\begin{rem}
The univariate cases $n=1$, $m=0$ and $n=0$, $m=1$ are the Nasrallah-Rahman
integral and an integral identity of Gasper (equations (6.4.1) and (4.11.4)
of \cite{GasperG/RahmanM:2004}); the general $m=0$ case is due to Gustafson
\cite{GustafsonRA:1992}.
\end{rem}

We also obtain a nontrivial transformation by taking $a\sim p^\alpha$ for
$0<\alpha<1$, say $a=(pq)^{1/2}$ for symmetry.

\begin{thm}\label{thm:triglim_IC3}
For any nonzero parameters $t_0,\dots,t_n$, $u_0,\dots,u_{n+m+1}$,
$v_0,\dots,v_m$, $q$ with
\[
0<|q|<1,\quad
\prod_{0\le r\le n} t_r = \prod_{0\le r\le n+m+1}u_r\prod_{0\le r\le m}v_r
\]
and any complex parameter $w$, we have the limit
\begin{align}
\lim_{p\to 0}&
\prod_{0\le r<s\le n} \theta((pq)^{-1/2}t_rt_s;q)
I^{(m)}_{BC_n}(\frac{t_0}{(pq)^{1/4}},\dots,\frac{t_n}{(pq)^{1/4}},
               \frac{(pq)^{1/4}}{u_0},\dots,\frac{(pq)^{1/4}}{u_{m+n+1}},
               \frac{(pq)^{3/4}}{v_0},\dots,\frac{(pq)^{3/4}}{v_m}
;p,q)\notag\\
&{}=
\frac{(q;q)^n}{n!}
\int_{C^n}
\frac{\theta(\prod_{0\le r\le n} t_r/w\prod_{1\le i\le n} z_i;q)
      \prod_{1\le i\le n} \theta(z_i/w;q)}
     {\prod_{0\le r\le n} \theta(t_r/w;q)}
\prod_{1\le i<j\le n} \Gammat((z_i/z_j)^{\pm 1};q)^{-1}\notag\\
&\phantom{{}= \frac{(q;q)^n}{n!} \int_{C^n}}
\prod_{1\le i\le n}
  \frac{\prod_{0\le r\le n} \Gammat(t_r/z_i;q)
        \prod_{0\le r\le m+n+1} \Gammat(z_i/u_r;q)}
       {\prod_{0\le r\le m} \Gammat(z_iv_r;q)}
  \frac{dz_i}{2\pi\sqrt{-1}z_i}
\notag
\end{align}
where the contour contains all points of the form $q^it_r$,
$i\ge 0$, and excludes all points of the form $q^{-i}u_r$, $i\ge 0$.
\end{thm}

\begin{cor}
The trigonometric integral of Theorem \ref{thm:triglim_IC3} is independent
of $w$, and if multiplied by
\[
\prod_{0\le r\le n} \prod_{0\le s\le n+m+1} \Gammat(t_r/u_s;q)^{-1},
\]
is invariant under the involution
\[
(m,n;\dots,t_r,\dots;\dots,u_r,\dots;\dots,v_r\dots)
\to
(n,m;\dots,v_r\dots;\dots,u_r^{-1},\dots;\dots,t_r,\dots).
\]
\end{cor}

\begin{rem}
This can also be obtained as a limit of Corollary \ref{cor:triglim_IC12}
after first breaking the symmetry of the left-hand side as in Lemma
\ref{lem:nonsym_ICn}.  In particular, this should perhaps be thought of as
a degeneration rather than a direct limit; we mention it to point out
that that distinction is somewhat artificial (any degeneration should be
obtainable as a limit directly from the elliptic level), but more
importantly because the Type II analogue has important consequences.
\end{rem}

For the Type II integral, we again have a trivial limit in one case.

\begin{thm}\label{thm:triglim_IIC1}
For any parameters $u_0,\dots,u_{m+4}$, $v_0,\dots,v_m$, $q$ satisfying
\[
|q|<1,\quad t^{2n-2}\prod_{0\le r\le m+4} u_r = \prod_{0\le r\le m} v_r,
\]
we have the limit
\begin{align}
\lim_{p\to 0}
\II^{(m)}_{BC_n}&(u_0,\dots,u_{m+4},pq/v_0,\dots,pq/v_m;t;p,q)\notag\\
&=
\frac{(q;q)^n\Gammat(t;q)^n}{2^n n!}
\int_{C^n}
\prod_{1\le i<j\le n} 
\frac{\Gammat(t z_i^{\pm 1} z_j^{\pm 1};q)}
     {\Gammat(  z_i^{\pm 1}z_j^{\pm 1};q)}
\prod_{1\le i\le n}
\frac{\prod_{0\le r\le m+4} \Gammat(u_r z_i^{\pm 1};q)}
     {\Gammat(z_i^{\pm 2};q) \prod_{0\le r\le m} \Gammat(v_r z_i^{\pm 1};q)}
\frac{dz_i}{2\pi\sqrt{-1}z_i},\notag
\end{align}
where the contour $C$ satisfies $C=C^{-1}$, and for all $i\ge 0$,
contains the points $q^i u_r$ as well as the contour $q^i t C$.
\end{thm}

We recall from \cite{Rains:Transformations} the following identity, which
plays the role of Lemma \ref{lem:diffop_ICn} for the type II $BC_n$
integral.
\[
\prod_{1\le i\le n} (1+R(z_i))
\frac{\theta(u_0 z_i,u_1 z_i,u_2 z_i,t^{n-1}u_0u_1u_2/z_i;q)}
     {\theta(z_i^2;q)}
\prod_{1\le i<j\le n} \frac{\theta(tz_iz_j;q)}
                           {\theta(z_iz_j;q)}
=
\prod_{0\le i<n} \theta(t^i u_0u_1,t^i u_0u_2,t^i u_1u_2;q).
\]

\begin{lem}\label{lem:nonsym_IICn}
For any nonzero parameters $t_0,t_1$, $u_0,\dots,u_{m+1}$,
$v_0,\dots,v_m$, $p$, $q$ with
\[
0<|p|,|q|<1,\quad
t^{2n-2}t_0t_1 = \prod_{0\le r\le m+2}u_r\prod_{0\le r\le m}v_r
\]
and any complex parameters $a,w\ne 0$, we have the
identity
\begin{align}
\prod_{0\le i<n} &\theta(t^i t_0t_1/a;q)
\II^{(m)}_{BC_n}(t_0/a^{1/2},t_1/a^{1/2},
                 a^{1/2}/u_0,\dots,a^{1/2}/u_{m+2},
                 pq/a^{1/2}v_0,\dots,pq/a^{1/2}v_m
;t;p,q)\notag\\
{}={}&
\frac{(p;p)^n(q;q)^n\Gammae(t;p,q)^n}{n!}
\int_{C^n}
\prod_{1\le i\le n} \frac{\theta(z_i/w,t^{n-1}t_0t_1/wz_i;q)}
                         {\theta(t^{i-1} t_0/w,t^{i-1} t_1/w;q)}
\prod_{1\le i<j\le n}
  \frac{\Gammae(tp z_iz_j/a,pq z_iz_j/ a,t (z_i/z_j)^{\pm 1};p,q)}
       {\Gammae( p z_iz_j/a,pq z_iz_j/ta,  (z_i/z_j)^{\pm 1};p,q)}\notag\\
&
\phantom{\frac{(p;p)^n(q;q)^n\Gammae(t;p,q)^n}{n!}\int_{C^n}}
\prod_{1\le i\le n}
  \theta(p z_i^2/a;p)
  \Gammae(p t_0 z_i/a,p t_1 z_i/a,t_0/z_i,t_1/z_i;p,q)\notag\\
&\phantom{\frac{(p;p)^n(q;q)^n\Gammae(t;p,q)^n}{n!}\int_{C^n}\prod_{1\le
      i\le n}}
  \prod_{0\le r\le m+2} \frac{\Gammae(z_i/u_r;p,q)}
                             {\Gammae(pq z_iu_r/a;p,q)}
  \prod_{0\le r\le m} \frac{\Gammae(pq z_i/a v_r;p,q)}
                           {\Gammae(z_iv_r;p,q)}
\frac{dz_i}{2\pi\sqrt{-1}z_i},
\notag
\end{align}
where the contour $C$ is chosen so that for all $i,j\ge 0$, it contains the
points and contours
\[
p^i q^j t_0,p^i q^j t_1,p^i q^j a/u_r,p^{i+1}q^{j+1}/v_r,
\quad p^i q^j t C,p^i q^j t a/C
\]
and excludes the points and contours
\[
a/p^{i+1}q^jt_0,a/p^{i+1}q^jt_1,u_r/p^iq^j,av_r/p^{i+1}q^{j+1},
\quad
a/p^{i+1}q^jtC,C/p^iq^jt.
\]
\end{lem}

\begin{rem}
It is possible to choose a contour of the given form satisfying
$C=aC^{-1}$, namely $a^{-1}C_0$ where $C_0$ is a suitable contour for the
left-hand side.
\end{rem}

\begin{thm}\label{thm:triglim_IIC2}
For any nonzero parameters $u_0,\dots,u_{2m+3}$, $v_0$, $v_1$, $q$ with
\[
0<|q|<1,\quad
\prod_{0\le r\le 2m+3}u_r = t^{2n-2}v_0v_1
\]
and any complex parameters $w\ne 0$, we have the
identity
\begin{align}
\lim_{p\to 0}
\prod_{0\le i<n} &\theta(t^i v_0v_1/pq;q)
\II^{(m)}_{BC_n}(v_0/(pq)^{1/2},v_1/(pq)^{1/2},
                 (pq)^{1/2}/u_0,\dots,(pq)^{1/2}/u_{2m+3}
;t;p,q)\notag\\
{}={}&
\frac{(q;q)^n\Gammat(t;q)^n}{n!}
\int_{C^n}
\prod_{1\le i\le n} \frac{\theta(z_i/w,t^{n-1}v_0v_1/wz_i;q)}
                         {\theta(t^{i-1} v_0/w,t^{i-1} v_1/w;q)}
\prod_{1\le i<j\le n}
  \frac{\Gammat(tz_iz_j/q,z_iz_j  ,t (z_i/z_j)^{\pm 1};q)}
       {\Gammat( z_iz_j/q,z_iz_j/t,  (z_i/z_j)^{\pm 1};q)}\notag\\
&
\phantom{\frac{(q;q)^n\Gammat(t;q)^n}{n!}\int_{C^n}}
\prod_{1\le i\le n}
  (1-z_i^2/q)
  \Gammat(v_0 z_i/q,v_1 z_i/q,v_0/z_i,v_1/z_i;q)
  \prod_{0\le r\le 2m+3} \frac{\Gammat(z_i/u_r;q)}
                              {\Gammat(z_i u_r;q)}
\frac{dz_i}{2\pi\sqrt{-1}z_i}
\notag
\end{align}
where the contour $C$ is chosen so that for all $i\ge 0$, it contains the
points and contours
\[
q^i v_0,q^i v_1, \quad q^i t C,
\]
and excludes the points and contours
\[
q^{1-i}/v_0,q^{1-i}/v_1,u_r/q^i, \quad q^{1-i}/tC,q^{-i}C/t,
\]
assuming such a contour exists.
\end{thm}

\begin{rem}
It is easy to verify that there exist choices of the parameters for which a
circular contour of radius $q^{1/2}$ satisfies the given conditions, and
thus the integral on the right has a well-defined meromorphic extension to
general parameters (and the limit will continue to hold); the only question
is whether this can be obtained from a domain of integration of the form
$C^n$.
\end{rem}

\begin{cor}
When $m=0$, the above trigonometric integral evaluates to
\[
\prod_{0\le i<n} 
\frac{\Gammat(t^{i+1},t^i v_0 v_1/q;q)
      \prod_{0\le r\le 3} \Gammat(t^i v_0/u_r,t^i v_1/u_r;q)}
     {\prod_{0\le r<s\le 3} \Gammat(t^{-i} u_ru_s;q)}
\]
\end{cor}

We also obtain a transformation.

\begin{cor}
When $m=1$, the trigonometric integral of Theorem \ref{thm:triglim_IIC1} is equal to
\[
\prod_{0\le i<n}
\frac{
  \prod_{0\le r<s\le 5} \Gamma(t^i u_ru_s;p,q)}
{\prod_{0\le r\le 5} \Gamma(v_0/t^i u_r,v_1/t^i u_r;p,q)
  \Gamma(v_0v_1/t^i q;p,q)}
\]
times the image of the trigonometric integral of Theorem
\ref{thm:triglim_IIC2} under the specialization $u_i\mapsto t^{(n-1)/2}
u_i$, $v_i\to t^{-(n-1)/2}v_i$.
\end{cor}

There are other transformations relating these integrals, but all can be
obtained by applying the above transformation to one or both sides of a
transformation of the integral of Theorem \ref{thm:triglim_IIC1} alone.

Similarly, taking $a=p^{1/2}$ above, we obtain the limit

\begin{thm}\label{thm:nonsym_IICn}
For any nonzero parameters $t_0,t_1$, $u_0,\dots,u_{m+1}$,
$v_0,\dots,v_m$, $p$, $q$ with
\[
0<|p|,|q|<1,\quad
t^{2n-2}t_0t_1 = \prod_{0\le r\le m+2}u_r\prod_{0\le r\le m}v_r
\]
and any complex parameter $w\ne 0$, we have the
identity
\begin{align}
\lim_{p\to 0}
\prod_{0\le i<n}\theta(t^i t_0t_1/p^{1/2};q)
\II^{(m)}_{BC_n}&(t_0/p^{1/4},t_1/p^{1/4},
                 p^{1/4}/u_0,\dots,p^{1/4}/u_{m+2},
                 p^{3/4}q/v_0,\dots,p^{3/4}q/v_m
;t;p,q)\notag\\
{}={}
\frac{(q;q)^n\Gammat(t;q)^n}{n!}
\int_{C^n}&
\prod_{1\le i\le n} \frac{\theta(z_i/w,t^{n-1}t_0t_1/wz_i;q)}
                         {\theta(t^{i-1} t_0/w,t^{i-1} t_1/w;q)}
\prod_{1\le i<j\le n}
  \frac{\Gammat(t (z_i/z_j)^{\pm 1};q)}
       {\Gammat(  (z_i/z_j)^{\pm 1};q)}\notag\\
&\prod_{1\le i\le n}
  \frac{
  \Gammat(t_0/z_i,t_1/z_i;q)
  \prod_{0\le r\le m+2} \Gammat(z_i/u_r;p,q)}
  {\prod_{0\le r\le m}\Gammat(z_iv_r;p,q)}
\frac{dz_i}{2\pi\sqrt{-1}z_i},
\notag
\end{align}
where the contour $C$ is chosen so that for all $j\ge 0$, it contains the
points and contours
\[
q^j t_0,q^j t_1,
\quad q^j t C
\]
and excludes the points and contours
\[
u_r/q^j, \quad C/q^jt.
\]
\end{thm}

Of course, with the Type II integral, we are particularly interested in the
effect of multiplying the integrand for $m=0$ by the biorthogonal
functions.  Note that since the integral is taken over a compact curve in
each case, the limiting relation will continue to hold as long as the limit
of biorthogonal functions exists, and (more difficult) the revised contour
conditions are satisfiable in the limit.  The primary constraint is that we
may only consider $p$-abelian biorthogonal functions, since otherwise the
contour must contain at least one point converging to $\infty$ as $p\to
0$.  For the first two limits, there is no difficulty with convergence of
the biorthogonal function.  Indeed, the $p$-abelian biorthogonal functions
satisfy the further identities
\begin{align}
\tilde{\cal R}^{(n)}_{0\lambda}&(\dots
p^{\pm 1/2}z_i\dots;p^{1/2}t_0{:}p^{1/2}t_1,p^{-1/2}t_2,p^{-1/2}t_3;p^{1/2}u_0,p^{-1/2}u_1;t;p,q)\notag\\
&{}=
\tilde{\cal R}^{(n)}_{0\lambda}(\dots
p^{\pm 1/2}z_i\dots;p^{-1/2}t_0{:}p^{-1/2}t_1,p^{1/2}t_2,p^{1/2}t_3;p^{1/2}u_0,p^{-1/2}u_1;t;p,q)\notag\\
&{}=
\tilde{\cal R}^{(n)}_{0\lambda}(\dots
z_i\dots;t_0{:}t_1,t_2,t_3;u_0,u_1;t;p,q)
\notag
\end{align}
and thus in each case the relevant limit of biorthogonal functions is the
same.

For the $a=\sqrt{p}$ limit, the situation is more delicate, but we find
that if $t^{2n-2}a_0a_1b_0b_1cd=q$, then we have a well-defined limit
\begin{align}
R^{(n)}_{\lambda;AS\text{-}I}&(\dots,z_i,\dots;a_0{:}a_1,b_0,b_1;c,d;q,t;p)
\notag\\
&:=
\lim_{p\to 0}
\tilde{\cal R}^{(n)}_{0\lambda}(\dots p^{-1/4}z_i\dots;
p^{-1/4}a_0{:}p^{-1/4}a_1,p^{1/4}b_0,p^{1/4}b_1;
p^{1/4} c,p^{3/4} d;t;p,q)\notag\\
&=
\lim_{p\to 0}
\tilde{\cal R}^{(n)}_{0\lambda}(\dots p^{-1/4}/z_i\dots;
p^{1/4}a_0{:}p^{1/4}a_1,p^{-1/4}b_0,p^{-1/4}b_1;
p^{3/4} c,p^{1/4} d;t;p,q).
\notag
\end{align}
This is a multivariate analogue of the biorthogonal rational functions of
Al-Salam and Ismail \cite{AlSalamWA/IsmailMEH:1994}.  More precisely,
by specializing the $a=\sqrt{p}$ limit appropriately (in particular,
$w=c^{-1}$), we find that the functions
\[
R^{(n)}_{\lambda;AS\text{-}I}(\dots,z_i,\dots;a_0{:}a_1,b_0,b_1;c,d;q,t;p)
\text{ and }
R^{(n)}_{\mu;AS\text{-}I}(\dots,1/z_i,\dots;b_0{:}b_1,a_0,a_1;d,c;q,t;p)
\]
are biorthogonal with respect to the density
\[
\prod_{1\le i<j\le n}
\frac{\Gammat(t (z_i/z_j)^{\pm 1};q)}{\Gammat((z_i/z_j)^{\pm 1};q)}
\prod_{1\le i\le n}
\frac{
\Gammat(a_0/z_i,a_1/z_i,b_0 z_i,b_1 z_i;q)
}{
\Gammat(q/cz_i,qz_i/d,t^{n-1}a_0a_1c/z_i,t^{n-1}b_0b_1dz_i;q)
},
\]
which becomes Al-Salam and Ismail's density when
$t^{n-1}a_0a_1c=q^{1/2}=t^{n-1}b_0b_1d$ and $n=1$.  The constraints on the
contour are independent of $\lambda$, $\mu$, $c$ and $d$, and are simply
that $C$ must contain the points $q^i a_r$ and the contours $q^i t C$, and
exclude the points $1/q^i b_r$ and the contours $C/q^i t$.
If we then take $a_1,b_1\to
0$ and set $a_0=b_0=q^{1/2}$, we obtain polynomials biorthogonal with
respect to the density
\[
\prod_{1\le i<j\le n}
\frac{\Gammat(t (z_i/z_j)^{\pm 1};q)}{\Gammat((z_i/z_j)^{\pm 1};q)};
\]
these are, of course, simply the ordinary Macdonald polynomials
\cite{MacdonaldIG:1995}, up to a suitable normalization.  That these arise
as limits of the biorthogonal functions is not particularly new (since they
are limits of Koornwinder polynomials); what {\em is} new is that a limit
exists that respects the inner product.

It should be possible to obtain similar limits in the $A_n$ case; since the
contour conditions are significantly more complicated in that case,
however, we mention only the identity which presumably plays the role of
Lemma \ref{lem:diffop_ICn} in this case:
\begin{align}
\text{symm}_{S_{n+1}}
\frac{\prod_{0\le i\le n}
        \theta(x z_i/\prod_{0\le j<i} v_j\prod_{i<j\le n}u_j;q)
        \prod_{0\le j<i} \theta(v_j z_i;q)
        \prod_{i<j\le n} \theta(u_j z_i;q)}
{\prod_{0\le i<j\le n} z_j\theta(z_i/z_j;q)u_j\theta(v_i/u_j;q)}
&\notag\\
{}=
\theta(x\prod_{0\le j\le n}z_j;q)
\prod_{1\le i\le n} &\theta(x/\prod_{0\le j<i} v_j\prod_{i\le j\le n}u_j;q),
\notag
\end{align}
a special case of Theorem 4.4 of \cite{RosengrenH/SchlosserM:2005}.

\section{Rational limits}

The rational limit is most naturally viewed as a combination of the
hyperbolic and trigonometric limits, and thus in particular requires both
the asymptotic calculations from the hyperbolic case and the symmetry
breaking from the trigonometric case.  In addition, it can be reached by
taking $\omega_2\to 0$ in the hyperbolic case, or $q\to 1$ in the
trigonometric or elliptic cases.  We consider the limit from the elliptic
level, as the other levels introduce no further complications.  In each
case, the integrand factors as a product of $q$-theta functions and
functions to which Corollary \ref{cor:ell_to_rat3} applies; the exponential
behaviour of the integrand comes only from the former.

\begin{thm}\label{thm:ratlim_IC1}
For $\Im(\omega)>0$, $\sum_{0\le r\le 2n+m+2}\mu_r = \sum_{0\le r\le m}\nu_r$,
\begin{align}
\lim_{v\to 0^+}&
e(-n(2n+3)/24v\omega)
\left((p;p)\sqrt{2\pi v\omega/\sqrt{-1}}\right)^{n(2n+3)}\notag\\
&I^{(m)}_{BC_n}(e(v\mu_0),\dots,e(v\mu_{2n+m+2}),
               pe(v(\omega-\nu_0)),\dots,pe(v(\omega-\nu_m));p,e(v\omega))\notag\\
&=
\frac{(\sqrt{2\pi}\omega/\sqrt{-1})^{-n}}{2^n n!}
\int_{C^n}
\prod_{1\le i<j\le n} \Gammar(\pm x_i\pm x_j;\omega)^{-1}
\prod_{1\le i\le n}
  \frac{\prod_{0\le r\le 2n+m+2} \Gammar(\mu_r\pm x_i;\omega)}
       {\Gammar(\pm 2x_i;\omega)
        \prod_{0\le r\le m} \Gammar(\nu_r\pm x_i;\omega)}
  dx_i,
\notag
\end{align}
where $C$ is a contour agreeing with $\R$ outside a compact set, and
separating the points of the form $\mu_r+j\omega$, $j\ge 0$ from the points
of the form $-\mu_r-j\omega$, $j\ge 0$.
\end{thm}

\begin{proof}
We consider the case in which the original contour is the unit circle;
deformed cases are analogous.  As in the hyperbolic case, we make the
change of variables $z_i=e(x_i)$ and integrate over $[-1/2,1/2]^n$.  If we
divide the integrand by
\[
\prod_{1\le i\le j\le n} \theta(e(x_i+x_j);e(v\omega))
\prod_{1\le i<j\le n} \theta(e(x_i-x_j);e(v\omega))
\prod_{1\le i\le n} \theta(e(v\omega/2+x_i);e(v\omega))^{-2n-2}
\]
the remaining factors of the integrand are controlled by Corollary
\ref{cor:ell_to_rat3} to have at worst polynomial growth in $v$, and thus
the exponential behavior of the theta functions dominates.  In particular,
up to polynomial factors, the integrand satisfies the bound
\[
O(e((-1/2v\omega)[
\sum_{1\le i\le n} 2(n+1)\vartheta(x_i)
-
\sum_{1\le i<j\le n} \vartheta(x_i\pm x_j)
-
\sum_{1\le i\le n} \vartheta(2x_i)
]))
\]
which decays exponentially unless $x_1$,\dots,$x_n=o(1)$.  We can thus
restrict the integral to $[-1/4,1/4]^n$ and rescale the variables by $v$.
The result then follows from Corollary \ref{cor:ell_to_rat2}.
\end{proof}

The other case is more complicated, in that the exponential contribution to
the asymptotics is not enough to properly localize the integral.

\begin{thm}\label{thm:ratlim_IC2}
For $\Im(\omega)>0$, $\sum_{0\le r\le 2n+m+2}\mu_r = \sum_{0\le r\le m}\nu_r$,
\begin{align}
&\lim_{v\to 0^+}
e(-(m^2+2m)/24v\omega)
\left((p;p)\sqrt{2\pi v\omega/\sqrt{-1}}\right)^{2m^2+3m}
\prod_{0\le r<s\le m} \theta(e(v(\nu_r+\nu_s-\omega))/p;e(v\omega))\notag\\
&\phantom{\lim_{v\to 0^+}}I^{(n)}_{BC_m}(\sqrt{p}e(v(\omega/2-\mu_0)),\dots,\sqrt{p}e(v(\omega/2-\mu_{2n+m+2})),
               \frac{e(v(\nu_0-\omega/2))}{\sqrt{p}},\dots,\frac{e(v(\nu_m-\omega/2))}{\sqrt{p}};p,e(v\omega))\notag\\
&\qquad=
\frac{(\sqrt{2\pi}\omega/\sqrt{-1})^{-m}}{m!}
\int_{C^m}
\frac{\theta_h(\sum_{0\le r\le m} \nu_r-w-\sum_{1\le i\le m} x_i;\omega)}
     {\prod_{1\le i\le m} \theta_h(x_i-w;\omega)^{-1}\prod_{0\le r\le m} \theta_h(\nu_r-w;\omega)}
\frac{\prod_{1\le i\le j\le m} ((x_i+x_j)/\omega-1)}
{\prod_{1\le i<j\le m} \Gammar(\pm (x_i-x_j);\omega)}\notag\\
&\phantom{\qquad=\frac{(\omega/\sqrt{-1})^{-m}}{m!}\int_{C^m}{}}
\prod_{1\le i\le m}
  \prod_{0\le r\le 2n+m+2} \frac{\Gammar(x_i-\mu_r;\omega)}
                                {\Gammar(x_i+\mu_r;\omega)}
  \prod_{0\le r\le m} \Gammar(\nu_r+x_i-\omega,\nu_r-x_i;\omega)
  dx_i,
\notag
\end{align}
where $C$ is a contour agreeing with $\R$ outside a compact set, and
separating the points of the form $\nu_r+j\omega$, $j\ge 0$ from the points
of the form $\mu_r-j\omega$, $-\nu_r-(j-1)\omega$, $j\ge 0$.
\end{thm}

\begin{proof}
We begin with the integral of Lemma \ref{lem:nonsym_ICn}, replacing the
extra parameter by $e(vw)$.  The exponential factor in the asymptotics of
the resulting elliptic integrand is
\[
e((
 \vartheta(\sum_{1\le i\le m} x_i)
+\sum_{1\le i<j\le m}\vartheta(x_i-x_j)
-m\sum_{1\le i\le m}\vartheta(x_i)
)/2v\omega),
\]
and thus the integrand is exponentially small unless the sequence
$x_1,\dots,x_m$ in $\R/\Z$ interlaces (or nearly interlaces) with the
all-zero sequence.  More precisely, if we split the integral into $2^m$
integrals based on the decomposition $\R/\Z =[-1/4,1/4]\cup[1/4,3/4]$, then
any piece with more than one $[1/4,3/4]$ factor contributes an
exponentially small amount.  For the pieces with exactly one $[1/4,3/4]$
factor, we find that upon rescaling the $[-1/4,1/4]$ variables, the
resulting integrand has order $O(v^2)$ and is integrable; thus those pieces
again contribute a negligible amount to the limit.  We may thus restrict
our attention to $[-1/4,1/4]^n$, or equivalently (up to $O(v)$), the
integral over $[-1/4v,1/4v]$ of the rational limit integrand.  The omitted
tails are again either exponentially small or have integral of order
$O(v^2)$, so the result follows.
\end{proof}

\begin{cor}
The rational integral of Theorem \ref{thm:ratlim_IC1} is equal to
\[
\prod_{0\le r<s\le 2n+m+2} \Gammar(\mu_r+\mu_s;\omega)
\prod_{\substack{0\le r\le 2n+m+2\\0\le s\le m}} \Gammar(\nu_s-\mu_r;\omega)^{-1}
\prod_{0\le r<s\le m} \Gammar(\nu_r+\nu_s-\omega;\omega)^{-1}
\]
times the rational integral of Theorem \ref{thm:ratlim_IC2}.
\end{cor}

For the Type II integral, again the first case is straightforward.

\begin{thm}
For any parameters $\mu_0,\dots,\mu_{m+4}$, $\nu_0,\dots,\nu_m$, $\omega$,
$\tau$ satisfying
\[
\Im(\omega),\Im(\tau)>0,
\quad
(2n-2)\tau+\prod_{0\le r\le m+4} \mu_r = \prod_{0\le r\le m} \nu_r,
\]
we have the limit
\begin{align}
\lim_{v\to 0^+}
&
e(-n/4v\omega)
\left((p;p)\sqrt{2\pi v\omega/\sqrt{-1}}\right)^{2n(2n-3)\tau/\omega+6n}
\notag\\
&
\II^{(m)}_{BC_n}(e(v\mu_0),\dots,e(v\mu_{m+4}),
pe(v(\omega-\nu_0)),\dots,pe(v(\omega-\nu_m));
e(v\tau);p,e(v\omega))\notag\\
&{}=
\frac{\Gamma_r(\tau;\omega)^n}{(\sqrt{2\pi}\omega/\sqrt{-1})^n 2^n n!}
\int_{C^n}
\prod_{1\le i<j\le n}
\frac{\Gamma_r(\tau\pm x_i\pm x_j;\omega)}{\Gamma_r(\pm x_i\pm x_j;\omega)}
\prod_{1\le i\le n}
\frac{\prod_{0\le r\le m+4} \Gamma_r(\mu_r\pm x_i;\omega)}
     {\Gamma_r(\pm 2x_i;\omega)\prod_{0\le r\le m} \Gamma_r(\nu_r\pm
       x_i;\omega)}
dx_i
\notag
\end{align}
where the contour $C$ agrees with $\R$ outside a compact set, satisfies
$C=-C$ and for all $i\ge 0$, contains the points $i\omega+\mu_r$ as well as
the contour $i\omega+\tau+C$.
\end{thm}

\begin{proof}
The exponential factor in the asymptotics of the integrand is
\[
O(e(\sum_i (\vartheta(2x_i)-4\vartheta(x_i))/2v\omega)),
\]
which is exponentially small unless $x_i\equiv 0$.  The limit follows as
above.
\end{proof}

The nonsymmetric Type II integral has even worse behavior than the
nonsymmetric Type I case, however.

\begin{thm}
Let $\mu_0,\dots,\mu_{2m+3}$, $\nu_0$, $\nu_1$, $\omega$, $\tau$, $p$ be
parameters such that $|p|<1$, $\Im(\omega),\Im(\tau)>0$, and
\[
(2n-2)\tau+\nu_0+\nu_1=\sum_r \mu_r,
\]
as well as the convergence condition $\Re(\tau/\omega)>-1/n$.
Then we have the limit
\begin{align}
\lim_{v\to 0}&
e(-n/4v\omega)
\left((p;p)\sqrt{2\pi v\omega/\sqrt{-1}}\right)^{2n(2n-3)\tau/\omega+6n}
\prod_{0\le i<n} \theta(e(v(i\tau+\nu_0+\nu_1-\omega)/p);q)\notag\\
&\II^{(m)}_{BC_n}(
\frac{e(v(\nu_0-\omega/2)}{\sqrt{p}},\frac{e(v(\nu_1-\omega/2)}{\sqrt{p}},
\sqrt{p}e(v(\omega/2-\mu_0)),\dots,\sqrt{p}e(v(\omega/2-\mu_{2m+3}));
e(v\tau);p,e(v\omega))\notag\\
&{}=
\frac{\Gamma_r(\tau;\omega)^n}{(\sqrt{2\pi}\omega/\sqrt{-1})^n n!}
\int_{C^n}
\prod_{1\le i<j\le n}
  \frac{\Gamma_r(\tau+x_i+x_j-\omega,x_i+x_j,\tau+x_i-x_j,\tau+x_j-x_i;\omega)}
  {\Gamma_r(x_i+x_j-\omega,x_i+x_j-\tau,x_i-x_j,x_j-x_i;\omega)}\notag\\
&\phantom{
{}=
\frac{\Gamma_r(\tau;\omega)^n}{(\sqrt{2\pi}\omega/\sqrt{-1})^n n!}
\int_{C^n}}
\prod_{1\le i\le n}
  (2x_i/\omega-1)
  \Gamma_r(\nu_0+x_i-\omega,\nu_1+x_i-\omega,\nu_0-x_i,\nu_1-x_i;\omega)
\notag\\
&\phantom{
{}=
\frac{\Gamma_r(\tau;\omega)^n}{(\sqrt{2\pi}\omega/\sqrt{-1})^n n!}
\int_{C^n}
\prod_{1\le i\le n}}
  \frac{\theta_h(x_i-w,(n-1)\tau+\nu_0+\nu_1-w-x_i;\omega)}
       {\theta_h((i-1)\tau+\nu_0-w,(i-1)\tau+\nu_1-w;\omega)}
  \prod_{0\le r\le 2m+3} \frac{\Gammar(x_i-\mu_r;\omega)}{\Gammar(x_i+\mu_r;\omega)}
  dx_i,
\notag
\end{align}
where the contour $C$ is chosen so that for all $i\ge 0$, it contains the
points and contours
\[
i\omega+\nu_0,i\omega+\nu_1, \quad i\omega+\tau+C,
\]
and excludes the points and contours
\[
(1-i)\omega-\nu_0,(1-i)\omega-\nu_1,\mu_r-i\omega, \quad
(1-i)\omega-\tau-C,C-i\omega-\tau,
\]
assuming such a contour exists.
\end{thm}

\begin{proof}
The exponential factors in the asymptotics of the elliptic integrand
actually cancel completely, with the result that the integrand has
polynomial asymptotics.  The result will follow from dominated convergence
if we can show that the rational integrand converges.

For the rational tails, we find that (assuming $C=\omega/2+\R$ for
simplicity) the integrand converges iff the integral
\[
\int_{\R^n}
\prod_{1\le i<j\le n} |x_i^2-x_j^2|^{2\Re(\tau/\omega)}
\prod_{1\le i\le n} |x_i+\sqrt{-1}\epsilon|^{-4(n-1)\Re(\tau/\omega)-3} dx_i
\]
converges for $\epsilon>0$; equivalently, via the change of variables $y_i
= x_i^2/\epsilon^2$, the theorem reduces to the convergence of
\[
\int_{[0,\infty)^n}
\prod_{1\le i<j\le n} |y_i-y_j|^{2\Re(\tau/\omega)}
\prod_{1\le i\le n} |1+y_i|^{-2(n-1)\Re(\tau/\omega)-2} dy_i.
\]
Up to linear fractional transformation, this is an instance of the Selberg
integral (\cite{SelbergA:1944}, stated as Corollary \ref{cor:Selberg}
below), and thus converges as long as $\Re(\tau/\omega)>-1/n$, as required.
\end{proof}

\section{Classical limits}

The final limit we consider is that corresponding to the usual beta
integral.  Although the beta integral itself is generally viewed as the
bottom level, this is in fact a somewhat misleading view, as the integrals
we obtain are in fact still elliptic (involving powers of theta functions).
For the $BC_n$ cases, a suitable change of variables exists that essentially
eliminates the dependence on $p$, but the corresponding change of variables
for the $A_n$ integral is much less obvious (if it exists at all).
Furthermore, even for the beta case, the classical transformation analogue
can only easily be reached by degenerating either the hyperbolic or
elliptic levels; the symmetry breaking of the trigonometric and rational
cases introduces unnecessary complications.

Since the integrand involves powers of theta functions, there are in
general some subtle issues involving choices of branch.  It will thus be
convenient to restrict to the case $p$ real, where the phases are easier to
control.  We have the following.

\begin{lem}
Choose $p$ and $z$ such that $-1<p<1$ and $|z|=1$.  Then the standard branch of
$\log\theta(z;p)$ satisfies
\begin{align}
\log\theta(z;p) &= \log(-z)/2 + \log|\theta(z;p)|,\notag\\
\log\theta(p^{1/2}z;p) &= \log|\theta(p^{1/2}z;p)|.\notag
\end{align}
\end{lem}

\begin{proof}
We have
\[
\log\theta(z;p) = \log(1-z) + \sum_{1\le i} \log(1-p^i z)+\log(1-p^i/z),
\]
taking the principal branch of the logarithm on the right-hand side.  Now,
\[
\log(1-p^i z)+\log(1-p^i/z) = 2\log|1-p^i z|,
\]
so it suffices to show that
\[
\log(1-z) = \log(-z)/2 + \log|1-z|,
\]
which follows from the observation $(1-z)/\sqrt{-z}=|1-z|$.
Similarly,
\[
\log\theta(p^{1/2}z;p) = \sum_{0\le i} 2\log|1-p^{i+1/2}z|.
\]
\end{proof}

We will thus assume $-1<p<1$ in the sequel; note, however, that the case of
more general $p$ can be obtained by replacing
\begin{align}
|\theta(z;p)|^{\kappa}&\mapsto (-z)^{\kappa/2}
\theta(z;p)^\kappa\notag\\
|\theta(p^{1/2}z;p)|^{\kappa}&\mapsto \theta(p^{1/2}z;p)^\kappa.\notag
\end{align}
In any event, we need all parameters to have absolute value $1$, $|p|^{1/2}$,
or $|p|$ within $1+o(1)$ for this to work.

For the Type I $BC_n$ integral, it is particularly natural to take
$2m+2$ parameters to have norm $|p|^{1/2}$, at which point both
sides of the transformation take the same form.  More general cases
could be considered, but appear to give rise to the same limiting
identities, so we will restrict our attention to the simplest case.

Given points $x,y,z$ on the unit circle
with $x,z$ distinct, $y\in [x,z]$ or equivalently $x\le y\le z$ indicates
that $y$ is on the closed counterclockwise arc from $x$ to $z$, and
similarly for open arcs.

\begin{thm}
Let $a_0,\dots,a_n$, $b_0,\dots,b_m$ be points on the unit circle with 
\[
1\le a_0<a_1<\cdots<a_n\le -1,
\]
let $\omega$ be a point in the upper half-plane, and let
$\alpha^{\pm}_0,\dots,\alpha^{\pm}_n$, $\beta^{\pm}_0,\dots,\beta^{\pm}_m$
be parameters such that $\Re(\alpha^{\pm}_r)>0$ and
\[
\sum_{0\le r\le n} \alpha^+_r+\alpha^-_r=\sum_{0\le r\le m} \beta^+_r+\beta^-_r.
\]
Then, writing $q=e(v\omega)$, $\alpha_r=\alpha^+_r+\alpha^-_r$,
$\beta_r=\beta^+_r+\beta^-_r$, we have
\begin{align}
\lim_{v\to 0^+}&
\frac{\Gammae(q^{\sum_r \alpha_r};p,q)}
{
\prod_{0\le r\le n} \Gammae(q^{\alpha_r};p,q)
\prod_{0\le r<s\le n} \Gammae(q^{\alpha^+_r+\alpha^+_s} a_r a_s,q^{\alpha^+_r+\alpha^-_s} a_r/a_s,q^{\alpha^-_r+\alpha^+_s} a_s/a_r,q^{\alpha^-_r+\alpha^-_s}/a_r/a_s;p,q)
}
\notag\\
&I^{(m)}_{BC_n}(\dots,q^{\alpha^+_r} a_r,q^{\alpha^-_r}/a_r,\dots,
               \dots,p^{1/2} q^{1/2-\beta^+_r}/b_r,p^{1/2} q^{1/2-\beta^-_r} b_r,\dots;p,q)\notag\\
&=
\prod_{0\le r<s\le n} |\theta(a_r a_s^{\pm 1};p)|^{1-\alpha_r-\alpha_s}
\frac{\Gamma(\sum_{0\le r\le n} \alpha_r)}{\prod_{0\le r\le n} \Gamma(\alpha_r)}
(2\pi(p;p)^2)^n\notag\\
&\phantom{{}={}}
\int_{z_i\in [a_{i-1},a_i]}
\prod_{1\le i<j\le n} |\theta(z_i z_j^{\pm 1};p)|
\prod_{1\le i\le n}
  \frac{\prod_{0\le r\le n} |\theta(a_r z_i^{\pm 1};p)|^{\alpha_r-1}}
       {\prod_{0\le r\le m} |\theta(p^{1/2} b_r z_i^{\pm 1};p)|^{\beta_r}}
  \frac{|\theta(z_i^2;p)|dz_i}{2\pi\sqrt{-1}z_i}.
\notag
\end{align}
\end{thm}

\begin{proof}
Using Lemma \ref{lem:gen_tri_ellC}, we find that the integral decays
exponentially outside the stated product of arcs (or images under the
hyperoctahedral group).  The result then follows by dominated convergence.
\end{proof}

If we define $\phi(z) = -\theta(z;p)^2/\theta(-z;p)^2$, then
\begin{align}
\phi(z)-\phi(w) &= \frac{\theta(-1;p)^2 z\theta(wz^{\pm
    1};p)}{\theta(-z;p)^2 \theta(-w;p)^2}\notag\\
\phi'(z) &= \frac{\theta(-1;p)^2 (p;p)^2 \theta(z^2;p)}{\theta(-z;p)^4},
\notag
\end{align}
so for $z$ in the arc $[1,-1]$,
\[
d\phi(z) = 2\pi (p;p)^2 \frac{\theta(-1;p)^2}{|\theta(-z;p)|^4}
\frac{|\theta(z^2;p)|dz}{2\pi\sqrt{-1}z}
\]
and
\begin{align}
|\theta(a z^{\pm 1};p)|^\kappa
&=
|\phi(z)-\phi(a)|^\kappa
\frac{|\theta(-z;p)|^{2\kappa}|\theta(-a;p)|^{2\kappa}}{|\theta(-1;p)|^{2\kappa}}
\notag\\
|\theta(p^{1/2} b z^{\pm 1};p)|^\kappa
&=
|\phi(z)-\phi(p^{1/2}b)|^\kappa
\frac{|\theta(-z;p)|^{2\kappa} |\theta(-p^{1/2}b;p)|^{2\kappa}}
     {|\theta(-1;p)|^{2\kappa} }
\notag
\end{align}
Consequently, we can make a change of variables in the resulting
transformation to obtain the following result of Dixon \cite{DixonAL:1905}.

\begin{cor}\label{cor:Dixon}
For any parameters $a_0,\dots,a_n$, 
$b_0,\dots,b_m$, $\alpha_0,\dots,\alpha_n$, $\beta_0,\dots,\beta_m$ such that
\[
-b_m<-b_{m-1}<\dots<-b_0<a_0<a_1<\dots<a_n,\quad
0<\Re(\alpha_r),\Re(\beta_r),\quad
\sum_{0\le r\le n}\alpha_r = \sum_{0\le r\le m}\beta_r
\]
we have the identity
\begin{align}
\frac{\prod_{0\le i<j\le n} |a_i-a_j|^{1-\alpha_i-\alpha_j}}
     {\prod_{0\le i\le n} \Gamma(\alpha_i)\prod_{0\le j\le m} |a_i+b_j|^{\alpha_i}}
\int_{x_i\in [a_{i-1},a_i]}&
\prod_{1\le i<j\le n} |x_i-x_j|
\prod_{1\le i\le n}
\frac{\prod_{0\le j\le n} |a_j-x_i|^{\alpha_j-1}}
     {\prod_{0\le j\le m} |x_i+b_j|^{\beta_j}}
dx_i\notag\\
{}=
\frac{\prod_{0\le i<j\le m} |b_i-b_j|^{1-\beta_i-\beta_j}}
     {\prod_{0\le i\le m} \Gamma(\beta_i)\prod_{0\le j\le n} |b_i+a_j|^{\beta_i}}
&
\int_{x_i\in [b_{i-1},b_i]}
\prod_{1\le i<j\le m} |x_i-x_j|
\prod_{1\le i\le m}
\frac{\prod_{0\le j\le m} |b_j-x_i|^{\beta_j-1}}
     {\prod_{0\le j\le n} |x_i+a_j|^{\alpha_j}}
dx_i.
\notag
\end{align}
\end{cor}

A similar argument gives the type II analogue.

\begin{thm}\label{thm:ell_selberg_lim}
Let $a_0,a_1$, $b_0,\dots,b_m$ be points on the unit circle with 
\[
1\le a_0<a_1\le -1,
\]
let $\omega$ be a point in the upper half-plane, and let
$\alpha^{\pm}_0,\alpha^{\pm}_1$, $\beta^{\pm}_0,\dots,\beta^{\pm}_m$,
$\tau$ be parameters such that
\[
\Re(\alpha^{\pm}_0),\Re(\alpha^{\pm}_1),\Re(\tau)>0,\quad
2(n-1)\tau+\alpha_0+\alpha_1=\sum_{0\le r\le m} \beta_r,
\]
where $\alpha_r=\alpha^+_r+\alpha^-_r$, $\beta_r=\beta^+_r+\beta^-_r$.
Writing $q=e(v\omega)$, we have
\begin{align}
\lim_{v\to 0^+}
&\prod_{0\le i<n}
  \frac{\Gammae(q^{2(n-1)\tau + \alpha_0+\alpha_1-i\tau};p,q)}
       {\Gammae(q^{(i+1)\tau},q^{i\tau+\alpha_0},q^{i\tau+\alpha_1},
                 q^{i\tau+\alpha^+_0+\alpha^+_1}a_0a_1,
                 q^{i\tau+\alpha^+_0+\alpha^-_1}a_0/a_1,
                 q^{i\tau+\alpha^-_0+\alpha^+_1}a_1/a_0,
                 q^{i\tau+\alpha^-_0+\alpha^-_1}/a_0a_1;p,q)}
\notag\\
&\II^{(m)}_{BC_n}(q^{\alpha^+_0} a_0,q^{\alpha^-_0}/a_0,q^{\alpha^+_1} a_1,q^{\alpha^-_1}/a_1,
  \dots,p^{1/2} q^{1/2-\beta^+_r}/b_r,p^{1/2} q^{1/2-\beta^-_r}
  b_r,\dots;q^\tau,p,q)
\notag\\
&{}=
|\theta(a_0 a_1^{\pm 1};p)|^{n-n(n-1)\tau-n\alpha_0-n\alpha_1}
\prod_{0\le i<n}
  \frac{\Gamma(2(n-1)\tau + \alpha_0+\alpha_1-i\tau)\Gamma(\tau)}
       {\Gamma((i+1)\tau)\Gamma(i\tau+\alpha_0)\Gamma(i\tau+\alpha_1)}
\notag\\
&\phantom{{}={}}\frac{(2\pi(p;p)^2)^n}{n!}
\int_{[a_0,a_1]^n}
\prod_{1\le i<j\le n} |\theta(z_i z_j^{\pm 1};p)|^{2\tau}
\prod_{1\le i\le n}
  \frac{|\theta(a_0 z_i^{\pm 1};p)|^{\alpha_0-1}
        |\theta(a_1 z_i^{\pm 1};p)|^{\alpha_1-1}}
       {\prod_{0\le j\le m} |\theta(p^{1/2} b_j z_i^{\pm 1};p)|^{\beta_j}}
  \frac{|\theta(z_i^2;p)|dz_i}{2\pi\sqrt{-1}z_i}
\notag
\end{align}
\end{thm}

\begin{cor}\label{cor:Selberg} \cite{SelbergA:1944}
For any real numbers $a_0$, $a_1$, $b$ with $-b<a_0<a_1$, and parameters
$\alpha_0$, $\alpha_1$, $\tau$ with positive real part,
\begin{align}
\frac{1}{n!}
\int_{[a_0,a_1]^n}&
\prod_{1\le i<j\le n} |x_i-x_j|^{2\tau}
\prod_{1\le i\le n}
  \frac{|a_0-x_i|^{\alpha_0-1}
        |a_1-x_i|^{\alpha_1-1}}
       {|b+x_i|^{2(n-1)\tau+\alpha_0+\alpha_1}}
  dx_i
\notag\\
&=
\frac{|a_0-a_1|^{n(n-1)\tau+n\alpha_0+n\alpha_1-n}}
     {|a_0+b|^{n(n-1)\tau+n\alpha_1}|a_1+b|^{n(n-1)\tau+n\alpha_0}}
\prod_{0\le i<n}
  \frac{\Gamma((i+1)\tau)\Gamma(i\tau+\alpha_0)\Gamma(i\tau+\alpha_1)}
       {\Gamma(2(n-1)\tau + \alpha_0+\alpha_1-i\tau)\Gamma(\tau)}
\notag
\end{align}
\end{cor}

\begin{rem}
In fact (as observed in \cite{SelbergA:1944}), the constraint that
$\Re(\tau)>0$ is too strict, as can be seen from the fact that the
right-hand side remains finite and positive as long as
\[
\Re(\tau)>-1/n,-\Re(\alpha_0)/(n-1),-\Re(\alpha_1)/(n-1).
\]
One can presumably weaken the conditions of Theorem
\ref{thm:ell_selberg_lim} correspondingly.
\end{rem}

\begin{cor}
For any real numbers $a_0$, $a_1$, $b_0$, $b_1$ with $-b_1<-b_0<a_0<a_1$,
and parameters $\alpha_0$, $\alpha_1$, $\beta_0$, $\beta_1$, $\tau$ with
positive real part such that
$\alpha_0+\alpha_1=\beta_0+\beta_1$, we have the transformation
\begin{align}
\int_{[a_0,a_1]^n}
\prod_{1\le i<j\le n} |x_i-x_j|^{2\tau}
\prod_{1\le i\le n}&
\frac{|a_0-x_i|^{\alpha_0-1}|a_1-x_i|^{\alpha_1-1}}
     {|b_0+x_i|^{(n-1)\tau+\beta_0}|b_1+x_i|^{(n-1)\tau+\beta_1}}
dx_i\notag\\
=
C&
\int_{[b_0,b_1]^n}
\prod_{1\le i<j\le n} |x_i-x_j|^{2\tau}
\prod_{1\le i\le n}
\frac{|b_0-x_i|^{\beta_0-1}|b_1-x_i|^{\beta_1-1}}
     {|a_0+x_i|^{(n-1)\tau+\alpha_0}|a_1+x_i|^{(n-1)\tau+\alpha_1}}
dx_i
\notag
\end{align}
where
\[
C
=
\prod_{0\le i<n}
\frac{
|a_0+b_0|^{\alpha_0-\beta_0}
|a_0+b_1|^{\alpha_0-\beta_1}
|a_1+b_0|^{\alpha_1-\beta_0}
|a_1+b_1|^{\alpha_1-\beta_1}
}
{
|a_0-a_1|^{1-(n-1)\tau-\alpha_0-\alpha_1}
|b_0-b_1|^{-1+(n-1)\tau+\beta_0+\beta_1}}
\frac{\Gamma(i\tau+\alpha_0)\Gamma(i\tau+\alpha_1)}
     {\Gamma(i\tau+\beta_0)\Gamma(i\tau+\beta_1)}.
\]
\end{cor}

The $A_n$ case is similar.

\begin{thm}
Let $a_0,\dots,a_n$, $b_0,\dots,b_m$, $Z$ be points on the unit circle with 
\[
a_0<\dots<a_n<a_{n+1}:=a_0,
\]
let $\omega$ be a point in the upper half-plane, and let
$\alpha^{\pm}_0,\dots,\alpha^{\pm}_n$, $\beta^{\pm}_0,\dots,\beta^{\pm}_m$
be parameters such that
\[
\Re(\alpha^{\pm}_0),\dots,\Re(\alpha^{\pm}_n)>0,\quad
\sum_{0\le r\le n}\alpha_r-\sum_{0\le r\le m} \beta_r,
\]
where $\alpha_r=\alpha^+_r+\alpha^-_r$, $\beta_r=\beta^+_r+\beta^-_r$.
Writing $q=e(v\omega)$, we have
\begin{align}
\lim_{v\to 0^+}&
\frac{\Gammae(q^{\sum_r\alpha_r};p,q)}
{
\prod_{0\le r,s\le n} \Gammae(q^{\alpha^-_r+\alpha^+_s} a_s/a_r;p,q)
\Gammae(q^{\sum_r \alpha^-_r} Z/a_0\cdots a_n;p,q)
\Gammae(q^{\sum_r \alpha^+_r} a_0\cdots a_n/Z;p,q)
}\notag\\
&
I^{(m)}_{A_n}(Z|
  \dots,q^{\alpha^-_r}/a_r,\dots,
  \dots,p^{1/2}q^{1/2-\beta^-_r}b_r,\dots;
  \dots,q^{\alpha^+_r}a_r,\dots,
  \dots,p^{1/2}q^{1/2-\beta^+_r}/b_r,\dots;
p,q)\notag\\
&{}=
\prod_{0\le r<s\le n} |\theta(a_s/a_r;p)|^{1-\alpha_r-\alpha_s}
|\theta(Z/a_0\cdots a_n;p)|^{1-\sum_r \alpha_r}
\frac{\Gamma(\sum_r\alpha_r)
}
{
\prod_{0\le r\le n} \Gamma(\alpha_r)
}\notag\\
&\phantom{{}={}}
\frac{(2\pi(p;p)^2)^n}{(n+1)!}
\int_{\substack{\prod_{0\le i\le n} z_i = Z\\ z_i\in [a_i,a_{i+1}]}}
  \prod_{0\le i<j\le n}
     |\theta(z_i/z_j;p)|
  \prod_{0\le i\le n}
     \frac{\prod_{0\le r\le n} |\theta(z_i/a_r;p)|^{\alpha_r-1}}
          {\prod_{0\le r\le m} |\theta(p^{1/2}b_r z_i;p)|^{\beta_r}}
  \prod_{0\le i<n} \frac{dz_i}{2\pi\sqrt{-1}z_i}
\notag
\end{align}
\end{thm}

\begin{rem}
The fact that the $n+m+2$ theta functions
\[
\prod_{0\le i\le n} \theta(z_i/a_r;p),\ 
\prod_{0\le i\le n} \theta(p^{1/2} b_r z_i;p)
\]
for fixed $\prod_{0\le i\le n}z_i$ span an $n+1$-dimensional space is
presumably relevant to finding an appropriate change of variables to
eliminate the theta functions from the integrand.  Clearly, though, the
resulting computations would not give a {\em trivial} derivation of the
above integral from a more traditional multivariate beta integral.
\end{rem}

%
%
%

\bibliographystyle{plain}

\end{document}